\documentclass[12pt]{article}
\usepackage{latexsym,amsmath,amsfonts,amsbsy,curves,
bezier,epic,color,eepic,bezier,amscd,makeidx}

\newcommand{\qed}{\hbox{\rule[-2pt]
{3pt}{6pt}}}

\newtheorem{theo}{Theorem}[section]

\newtheorem{pro}[theo]{Proposition}
\newtheorem{dfe}[theo]{Definition}
 
\newtheorem{rem}[theo]{Remark}
\makeatletter

\@addtoreset{equation}{section}
\makeatother



\begin{title}
{On the existence of tight relative $2$-designs on binary Hamming association schemes}
\end{title}
\author{Eiichi Bannai,
Etsuko Bannai and
Hideo Bannai
}

\date{}
\begin{document}
\maketitle
\begin{abstract}
It is known that there is a close analogy between ``Euclidean $t$-designs vs. spherical 
$t$-designs" and ``Relative $t$-designs in binary Hamming association schemes vs. 
combinatorial $t$-designs". 
In this paper, we want to prove how much we can develop a similar theory in 
the latter situation, imitating the theory in the former one. We first prove that the 
weight function is 
constant on each shell for tight relative $t$-designs on $p$ shells on a wide class of 
Q-polynomial association schemes, including Hamming association schemes. 
In the theory of Euclidean $t$-designs on $2$ concentric spheres (shells), it is known 
that the structure of coherent configurations is naturally attached. However, it seems 
difficult to prove this claim in a general context. 
In the case of tight $2$-designs 
in combinatorial $2$-designs, there are great many tight $2$-designs, i.e., 
symmetric $2$-designs, 
while there are very few tight $2e$-designs for $e\geq 2$.
So, as a starting point,
we concentrate our study to the existence problem of tight relative $2$-designs, 
in particular on 
2 shells, in binary Hamming association schemes $H(n,2)$. 
We prove that every tight relative $2$-designs on 
2 shells in $H(n,2)$ has
the structure of coherent configuration. We determined all the possible
parameters of coherent configurations
attached to such tight relative $2$-designs for 
$n\leq 30$. Moreover
for each of them we determined whether 
there exists such a tight relative $2$-design or not, either by constructing them 
from symmetric $2$-designs 
or Hadamard matrices, or theoretically showing the non-existence. In particular, 
we show that 
for $n\equiv 6~ (\mbox{mod}~ 8)$, there exist such tight relative $2$-designs 
whose weight functions are not 
constant. These are the first examples of those with non-constant weight. 
\end{abstract}

 
 design, regular semi lattice.


\section{Introduction}
As is well known, there is a close analogy between the theory of combinatorial 
$t$-designs 
($t$-$(v,k,\lambda)$ designs) and the theory of spherical $t$-designs. 
Furthermore, it is known 
that there is a close analogy between the theory of Euclidean $t$-designs and 
the theory 
of relative $t$-designs in binary Hamming association schemes $H(n,2)$. 
Although this 
last analogy is known, it is not very well known up to now. 
(See, Delsarte \cite{Delsarte-1973, Delsarte-1977}, 
Delsarte-Seidel\cite{Delsarte-S-1989}, Bannai-Bannai \cite{Bannai-B-2012}.) 
The purpose of the present paper is to dig into more on this 
analogy. The theory of spherical harmonics has been developed into a very 
elaborate stage, and the 
theory is extremely beautiful. On the other hand, the theory of spherical 
functions on 
(Q-polynomial) association schemes is also developed, but in a sense it is 
more sophisticated. 
So, we need more careful treatments in order to get similar results known 
in Euclidean 
$t$-designs, for relative $t$-designs
on Q-polynomial, say, binary Hamming association schemes $H(n,2)$. 

For example, the tight spherical $t$-designs as well as the theory of tight 
combinatorial 
$t$-designs are both well studied, although complete classifications are not 
yet obtained 
at this stage. (Cf. \cite{Delsarte-G-S-1975}, \cite{Bannai-D-1}, 
\cite{Bannai-D-2}, \cite{Bannai-M-V-2004}, \cite{Nebe-V}, etc. for 
tight spherical $t$-designs, and 
Ray-Chaudhuri and Wilson \cite{Ray-Chaudhuri-W-1975}, 
Enomoto-Ito-Noda \cite{Enomoto-I-N-1979},
Bannai \cite{Bannai-1977}, Dukes-ShortGershman \cite{Dukes-S-2012}, etc. 
for 
tight combinatorial designs.) For tight Euclidean $t$-designs, the theory was 
developed in certain 
cases (\cite{Bannai-B-2006,Bannai-B-2009,Bannai-B-survey-2009,Bannai-B-2010,
EtBannai-2006, EtBannai-2009}). On the other hand, the theory of tight relative 
$t$-designs in $H(n,2)$ is less developed, so far. 
So, we will try to see how much the methods in Euclidean $t$-designs can be 
applied in the 
study of tight relative $t$-designs in $H(n,2).$  Here, we are mostly interested 
in the case where 
the number of spheres supporting the Euclidean designs or the number of shells 
supporting 
the combinatorial $t$-designs are relatively small,
equal to 2 in most cases. 
Also, we must put some strong restrictions 
on $t$, in some cases. 
In this paper, we want to obtain the following explicit results.\\
(i) We prove that for tight relative $2e$-designs in a Q-polynomial association 
scheme, the
weight function must be constant on each shell of the design, with a mild 
additional assumption. 
(See Theorem \ref{theo:2.1}). 
In the latter part of this paper, we restrict our study to tight relative 
$2e$-designs in $H(n,2)$, 
and also to $e=1.$ These are very strong restrictions, but still there are 
interesting examples 
and interesting theories. 
(ii) Using general theory of the study of tight relative $2e$-designs in 
$H(n,2)$, we determine all
the possible parameters of the tight relative $2$-designs on 
$2$-shells in $H(n,2),$ 
for explicit small values of $n,$ say $n\leq 30.$ Then, 
(iii) we determine the existence and the non-existence with those parameters 
listed 
in (ii). Very interesting feature is that we did find some examples of tight relative 
$2$-designs 
in $H(n,2),$ where the weight functions are not constant. It seems such 
examples were 
not known explicitly before. Here, we use some results obtained in 
\cite{Li-B-B-2013}. 

Our results obtained in this paper are only for special cases, but we expect that this 
approach will shed some light on the future studies of more general theory of 
(tight) relative $2e$-designs for bigger $e$ in more general Q-polynomial 
association schemes. 
As for the information on association schemes, e.g., more general P-polynomial, 
Q-polynomial
or P- and Q-polynomial schemes refer 
\cite{Bannai-I-1984}.

Now we introduce notation we use in this paper and some important definitions.
Let $\mathfrak X=(X,\{R\}_{0\leq r\leq d})$ be a symmetric association scheme.
Let $u_0\in X$ fixed arbitrarily. 
Let $X_r=\{x\in X\mid  (u_0,x)\in R_{r}\}$ for $r=0,1,\ldots,d$.
$X_0,X_1,\ldots,X_d$ are called shells of $\mathfrak X$.
$\mathcal F(X)$ be the vector space consists of all the real valued
functions on $X$. In the following argument we often identify 
$\mathcal F(X)$ with the vector space $\mathbb R^{|X|}$ indexed by $X$.
When we consider spherical designs or Euclidean designs,
we use the properties of vector spaces of polynomials. 
For the usual 
polynomials in $n$ variables defined on $\mathbb R^n$, it is
convenient to consider the subspaces $Hom_j(\mathbb R^n)$ spanned by 
all the homogeneous polynomials of degree $j$.
If $\mathfrak X$ is a P-polynomial scheme, it is natural to
consider the following subspace of $\mathcal F(X)$.
For any $z\in X_j$, we define $f_z\in \mathcal F(X)$ by
\begin{eqnarray}
f_z(x)=\left\{
\begin{array}{ll}
1&\mbox{if $x\in X_i$, $i\geq j$ and $(x,z)\in R_{i-j}$,}\\
0&\mbox{other wise}.
\end{array}
\right.
\end{eqnarray} 
Let $Hom_j(X)=\langle f_z\mid z\in X_j\rangle$. Then we have the following 
decomposition of $\mathcal F(X)$ into direct sum of subspaces.
\begin{eqnarray}
\mathcal F(X)=\mbox{Hom}_0(X)+\mbox{Hom}_1(X)+\cdots+\mbox{Hom}_d(X).
\label{equ:1.2}
\end{eqnarray}
Clearly we have
$$\dim(\mbox{Hom}_j(X))=|X_j|=k_j~(0\leq j\leq d).$$
When $\mathfrak X$ is a Q-polynomial scheme, it is natural to
consider the following subspace of $\mathcal F(X)$.
Let $E_0,E_1,\ldots,E_d$ be the primitive idempotents which give the
Q-polynomial structure of $\mathfrak X$.
For each $E_j$, let $L_j(X)$ be the subspace of $\mathcal F(X)$ 
spanned by all the column vectors of $E_j$.
Then we have $\dim (L_j(X))=rank(E_j)=m_j$ and
we have the following decomposition of $\mathcal F(X)$ into orthogonal
sum of subspaces. 
\begin{eqnarray}
\mathcal F(X)=L_0(X)\perp L_1(X)\perp \cdots\perp L_d(X).
\label{equ:1.3}
\end{eqnarray}
For each of the decomposition of $\mathcal F(X)$ given above we can
develop theory of relative $t$-designs for weighted subset 
$(Y,w)$ of $X$ using the similar setting as for the Euclidean designs. 
We use the following notation.
Let $\{r_1,r_2,\ldots,r_p\}=\{r\mid X_r\cap Y\not=\emptyset\}$.
Let $S=X_{r_1}\cup X_{r_2}\cup\cdots\cup X_{r_p}$.
We say $(Y,w)$ is supported by $p$ shells.
Let $Y_{r_i}= X_{r_i}\cap Y$, $i=1,2,\ldots, p$.
We also define 
$A(Y_{r_i},Y_{r_j})=\{\alpha\mid (x,y)\in R_\alpha, x\in Y_{r_i}, y\in Y_{r_j}, 
x\not=y\}$
for $1\leq i,j\leq p$.
We also use the notation $A(Y_{r_i})
=A(Y_{r_i},Y_{r_i})$ for $1\leq i\leq p$. 

As for the decomposition given by (\ref{equ:1.2}) for $H(n,2)$,
Delsarte-Seidel \cite{Delsarte-S-1989} defined the design as regular $t$-wise
balanced design. 

In this paper we consider the decomposition (\ref{equ:1.3}) for Q-polynomial 
schemes.
The concept of 
relative $t$-design with respect a fixed point $u_0\in X$ is related to the 
decomposition given by (\ref{equ:1.3}) for Q-polynomial schemes.
It was first defined by Delsarte in 1977 \cite{Delsarte-1977}. 
Without noticing his paper, we gave a definition of
relative $t$-designs with respect to $u_0\in X$
analyzing the concept of Euclidean $t$-designs
\cite{Bannai-B-2012}.
Later H. Tanaka informed us the existence of the paper
by Delsarte in 1977 \cite{Delsarte-1977}. In \cite{Bannai-B-2012}, 
we prove that our definition 
is equivalent to that of Delsarte.
We found that the theory of relative $t$-designs 
with respect to a fixed point is very similar to the concept of 
Euclidean design, in which the origin $0\in \mathbb R^n$
is a special point.

The following is the definition of relative $t$-design in the style of
Euclidean $t$-design (see \cite{Bannai-B-2012}).
\begin{dfe}
\label{def:1.1} Let $(Y,w)$ be a weighted subset of $X$ with positive 
weight function $w$ on $Y$.
$(Y,w)$ is called a relative $t$-design with respect to $u_0$ if the following 
condition holds.
\begin{eqnarray}
\sum_{i=1}^p\sum_{x\in X_{r_i}}
\frac{W_{r_i}}{|X_{r_i}|}f(x)
=\sum_{y\in Y}w(y)f(y)
\label{equ:1.4}
\end{eqnarray}
for any function $f\in L_0(X)\perp L_1(X)\perp
\cdots\perp L_{t}(X)$, where $W_{r_i}=\sum_{y\in Y_{r_i}}w(y)$, 
$i=1,2,\ldots,p$.
\end{dfe}
The following theorem is known \cite{Bannai-B-2012}.
\begin{theo}
Let $(Y,w)$ be a relative $2e$-design of a Q-polynomial scheme.
Then the following inequality holds.
\begin{eqnarray}
|Y|\geq \dim(L_0(S)+L_1(S)+\cdots+L_e(S)),
\label{equ:1.5}
\end{eqnarray}
where $L_j(S)=\{f|_S\mid f\in L_j(X)\}$, $j=0,1,\ldots, e$. 
\label{theo:1.2}
\end{theo}
\begin{dfe}
[tight relative $2e$-design with respect to $u_0$]
If equality holds in (\ref{equ:1.5}) in Theorem \ref{theo:1.2},
then $(Y,w)$ is called a tight relative $2e$-design with respect to $u_0$.
\end{dfe}
In the following we only consider
the nontrivial tight $2e$-designs $(Y,w)$.
That is, $Y$ does not contain $X_r$ for any $r$, $0\leq r\leq d$. 
\begin{rem} \begin{enumerate}
\item It is conjectured that 
$$\dim(L_0(S)+L_1(S)+\cdots+L_e(S))=m_e
+m_{e-1}+\cdots+m_{e-p+1}$$ 
holds for Q-polynomial schemes
with some trivial exceptions. For binary hamming scheme
it is proved that the conjecture is true 
\cite{Xiang-2012}.
\item In \cite{Delsarte-S-1989}, it is proved that a regular $2e$-wise 
balanced design
$(Y,w)$ satisfies 
\begin{eqnarray}
|Y|\geq \dim(\mbox{Hom~}_0(S)+\mbox{Hom~}_1(S)+\cdots
+\mbox{Hom~}_e(S)).
\label{equ:1.6}
\end{eqnarray}
However Delsarte-Seidel \cite{Delsarte-S-1989} mentioned that the
explicit computation of 
$\dim(\mbox{Hom~}_0(S)+\mbox{Hom~}_1(S)+\cdots+\mbox{Hom~}_e(S))$ will be difficult. Recently Xiang \cite{Xiang-2012} proved that
$$\dim(\mbox{Hom~}_0(S)+\mbox{Hom~}_1(S)+\cdots+\mbox{Hom~}_e(S))=k_e+k_{e-1}+\cdots+k_{e-p+1}$$
holds for $H(n,2)$. It is also proved that
$$\mbox{Hom~}_0(X)+\mbox{Hom~}_1(X)+\cdots+\mbox{Hom~}_e(X)
=L_0(X)+L_1(X)+\cdots+L_e(X)$$
for some P-and Q-polynomial schemes including $H(n,2)$ \cite{Bannai-B-S-T}.
Hence conjecture is correct for $H(n,2)$, if $S=X_{r_1}\cup X_{r_2}\cdots \cup X_{r_p}$ satisfies some suitable condition
to avoid the cases which trivially do not satisfy the conjecture.
\end{enumerate}
\end{rem}
In \S 2, we give our main results. In \S3 and \S4, we give the proofs of the main results.

\section{Main theorems}
\begin{theo}
Let $\mathfrak X=(X,\{R_r\}_{0\leq r\leq d})$
be a Q-polynomial scheme.
Let $G$ be the automorphism group of $\mathfrak X$.
Let $(Y,w)$ be a tight relative $2e$-design with respect to $u_0$
supported by $p$ shells.
Assume that the stabilizer $G_{u_0}$ of $u_0$
acts transitively on every shell $X_r,~ 1\leq r\leq d$.
Then the weight function $w$ of any tight relative $2e$-design 
$(Y,w)$ is constant on each 
$Y_{r_i}~(1\leq i\leq  p)$.
\label{theo:2.1}
\end{theo}
\begin{theo}
Let $(Y,w)$ be a tight relative $2$-design of the binary
Hamming scheme $H(n,2)$ 
supported by 2 shells, $S=X_{r_1}\cup X_{r_2}$. 
Let $N_{r_i}=|Y_{r_i}|$, $w(y)=w_{r_i}$ on 
$y\in Y_{r_i}$ for $i=1,2$.  
\begin{enumerate}
\item For any integers $r_1,~r_2$ satisfying
$1\leq r_1<r_2\leq n-1$, the following holds 
$$|A(Y_{r_1})|=|A(Y_{r_2})|
=|A(Y_{r_1},Y_{r_2})|=1.$$ 
This means that $Y=Y_{r_1}\cup Y_{r_2}$ has a structure of coherent configuration.
\item Assume
$1\leq r_1<r_2\leq n-1$
and $n\leq 30$, then 
the set of parameters $\{n,r_1,r_2,N_{r_1},\\
N_{r_2},
\alpha_1,
\alpha_2, \gamma,\frac{w_{r_2}}{w_{r_1}}\}$
is among those listed in {\rm \S 4.4},
here $A(Y_{r_i})=\{\alpha_i\}$ for $i=1,2$ and $A(Y_{r_1},Y_{r_2})=\{\gamma\}$.
\item If $n\equiv 6~(\mbox{mod}~8)$, and there exists
Hadamard matrix of size $\frac{1}{2}n+1$, then
there exists tight relative 2 design 
$Y\subset X_2\cup X_{\frac{n}{2}}$
($r_1=2$, $r_2=\frac{n}{2}$) whose
weights satisfy 
$\frac{w_{\frac{1}{2}n}}{w_2}=\frac{8}{n+2}$, that is, 
$w$ is not constant on $Y$. 
\item If $n\leq 30$, $w_{r_1}\neq w_{r_2}$ and
$n\not\equiv 6 (\mbox{mod}~ 8)$
or $n\equiv 6 (\mbox{mod}~ 8)$ and $Y$ is not related to the Hadamard matrices in (3), then there is no tight relative $2$-designs with respect to $u_0$.

\end{enumerate}

\label{theo:2.2}
\end{theo}
\begin{rem}
Since we consider only nontrivial tight 
designs and $|X_0|=|X_n|=1$, we may assume $1\leq r_1<r_2\leq n-1$. 
\end{rem}
In \S 3, we give the proof of Theorem \ref{theo:2.1}.
In \S 4, we give the proof of Theorem \ref{theo:2.2}.
Proposition \ref{pro:4.3} in \S 4.1 gives the explicit 
formula for $\frac{w_{r_2}}{w_{r_2}}$, $\alpha_1, \alpha_2$, and
$\gamma$ in terms of 
$n$, $r_1$, $r_2$, $N_{r_1}$. This implies Theorem \ref{theo:2.2} (1).
Proposition \ref{pro:4.1} and Proposition \ref{pro:4.2} in \S 4.1 give 
some useful formulas to prove Proposition \ref{pro:4.3}.
In \S 4.2 we give the proof of Propositions \ref{pro:4.1}, \ref{pro:4.2} and
\ref{pro:4.3}. 
In \S 4.4 we give the table of possible sets of parameters
$\{n, r_1, r_2,N_{r_1},
N_{r_2},  \alpha_1,
\alpha_2,  \gamma,\frac{w_{r_2}}{w_{r_1}}\}$ for tight relative 2-design $(Y,w)$
with respect to $u_0$ for $n\leq 30$. We give
two kind of construction theorems.
One is the construction from Hadamard matrices 
and the other is the construction from symmetric designs (Proposition \ref{pro:4.5} in \S 4.4).

To obtain the feasible parameters,
$\{n, r_1, r_2,N_{r_1},
N_{r_2},  \alpha_1,
\alpha_2,  \gamma,\frac{w_{r_2}}{w_{r_1}}\}$ 
in the table given in \S 4, we mainly used the properties of Q-polynomial 
structure of $H(n,2)$.
In \cite{Delsarte-1977}, Delsarte proved that if association scheme is attached to a regular semi lattice (then it is a P-polynomial scheme), and if it
also has Q-polynomial structure with it's ordering, then $(Y,w)$ is a
relative $t$-design with respect to $u_0$
if and only if it is a geometric relative $t$-design
with respect to $u_0$.
For $H(n,2)$, a geometric relative $t$-design
with respect to $u_0$ is nothing but
a regular $t$-wise balanced design. In \S 4.3 we briefly introduce 
regular semi-lattices and geometric relative $t$-designs. 
We also use the property of regular $t$-wise balanced design 
(Proposition \ref{pro:4.4} in \S4) to show the non-existence 
of such a design for some feasible parameters in the
table.
\begin{rem}
We conclude this section by mentioning that Woodall \cite{Woodall-1970},
in particular Theorem 8 in \cite{Woodall-1970} essentially discuss similar 
problem as ours under the additional assumption that the weight function
is constant. It would be interesting to compare our approach with that of 
Woodall \cite{Woodall-1970}    
\end{rem}
 
\section{Proof of Theorem \ref{theo:2.1}}
Let $L(S)$ be the vector space of real valued 
functions on $S$. 
We consider the inner product on $L(S)$ 
defined for $f,g\in L(S)$ by
\begin{eqnarray}
\langle f,g\rangle
=\sum_{i=1}^p\frac{W_{\nu_i}}{|X_{\nu_i}|}
\sum_{x\in X_{\nu_i}}
f(x)g(x).
\label{equ:3.1}
\end{eqnarray}
Let $\{\varphi_1,\ldots,\varphi_N\}\subset 
L_0(X)\perp L_1(X)\perp
\cdots\perp L_{e}(X)$.
Assume that 
$\{\varphi_1|_S,\ldots,\varphi_N|_S\}$
is an orthonormal basis of 
$L_0(S)+ L_1(S)+\cdots+ L_e(S)$
with respect to this inner product.
Let $H$ be the matrix whose rows are 
indexed by $Y$ with $N$ columns whose $(y,i)$-entry
is defined by $\sqrt{w(y)}\varphi_i(y)$. 
Since 
$fg\in  
L_0(X)\perp L_1(X)\perp\cdots\perp L_{2e}(X)$ 
holds 
for any $f,g\in L_0(X)\perp L_1(X)\perp\cdots
\perp L_e(X)$,
$\varphi_i\varphi_j\in L_0(X)\perp L_1(X)\perp\cdots
\perp L_{2e}(X)$.
Then we have the following
\begin{eqnarray}&&
(^tH~H)(i,j)=\sum_{y\in Y}w(y)\varphi_i(y)
\varphi_j(y)\nonumber\\
&&=\sum_{i=1}^p\sum_{x\in X_{r_i}}
\frac{W_{r_i}}{|X_{r_i}|}\varphi_i(x)
\varphi_j(x)=\delta_{i,j}.
\end{eqnarray}
This implies $\mbox{rank}(H)=|Y|\geq N=
\dim(L_0(S)+ L_1(S)+\cdots+ L_e(S)$.
If $|Y|=N$, then $H$ is an invertible 
matrix and $H~^tH=I$ holds.
Then we have
\begin{eqnarray}
(H^tH)(y_1,y_2)
=\sum_{i=1}^N\sqrt{w(y_1)w(y_2)}
\varphi_i(y_1)\varphi_i(y_2)=\delta_{y_1,y_2}.
\end{eqnarray}
This implies
\begin{eqnarray}
&& 
\sum_{i=1}^N\varphi_i(x)\varphi_i(y)=\delta(x,y)\frac{1}{w(y)}.
\end{eqnarray}
We introduce the following notation.
Let $\phi^{(j)}_u\in L_j(X)$ whose $x$-entry is
defined by $\phi^{(j)}_u(x)=\frac{1}{|X|}E_j(x,u)$ 
for $0\leq j\leq d$. Let 
$$\{\phi^{(j_1)}_{u_1},\phi^{(j_2)}_{u_2},\ldots,\phi^{(j_N)}_{u_N}\}\subset 
\bigcup_{j=0}^e\{\phi^{(j)}_u\mid u\in X\} $$
be a set of functions whose restrictions to $S$
forms a basis of $L_0(S)+ L_1(S)+\cdots+ L_e(S)$.
Let $u_s\in X_{l_s}$.
For simplicity let us write $\phi_s=\phi^{(j_s)}_{u_s}$
for $s=1,2,\ldots,N$.
From $\{\phi_1,\phi_2,\ldots,\phi_N\}$
we construct a set of $\{\varphi_1,\ldots,\varphi_N\}$ 
whose restrictions $\varphi_1|_S,\ldots,\varphi_N|_S$ 
to $S$ forms an orthonormal system in $L(S)$.
It is well known that Gram-Schmidt's method gives the
following formula for
$\varphi_1, \ldots,\varphi_N$.
\begin{eqnarray}
\varphi_j=\frac{1}{\sqrt{D_{j-1}D_j}}
\left|\begin{array}{cccc}
\langle\phi_1,\phi_1\rangle
&\langle\phi_2,\phi_1\rangle&\cdots
&\langle\phi_j,\phi_1\rangle\\
\langle\phi_1,\phi_2\rangle
&\langle\phi_2,\phi_2\rangle&\cdots
&\langle\phi_{j},\phi_2\rangle\\
\vdots&\cdots&\cdots&\vdots\\
\langle\phi_1,\phi_{j-1}\rangle
&\langle\phi_2,\phi_{j-1}\rangle&\cdots
&\langle\phi_{j},\phi_{j-1}\rangle\\
\phi_1&\phi_2&\cdots&\phi_j
\end{array}
\right|,\label{equ:3.5N}
\end{eqnarray}
where $D_j$ is the Gram determinant
given by
\begin{eqnarray}
D_j=
\left|\begin{array}{cccc}
\langle\phi_1,\phi_1\rangle
&\langle\phi_2,\phi_1\rangle&\cdots
&\langle\phi_j,\phi_1\rangle\\
\langle\phi_1,\phi_2\rangle
&\langle\phi_2,\phi_2\rangle&\cdots
&\langle\phi_{j},\phi_2\rangle\\
\vdots&\cdots&\cdots&\vdots\\
\langle\phi_1,\phi_{j-1}\rangle
&\langle\phi_2,\phi_{j-1}\rangle&\cdots
&\langle\phi_{j},\phi_{j-1}\rangle\\
\langle\phi_1,\phi_{j}\rangle
&\langle\phi_2,\phi_{j}\rangle&\cdots
&\langle\phi_{j},\phi_{j}\rangle
\end{array}
\right|.
\label{equ:3.6N}
\end{eqnarray}
The formula (\ref{equ:3.5N}) means
$\varphi_j$ is given by the linear sum of
$\phi_l$ with coefficient given by the
$(j,l)$-cofactor of the matrix given in 
(\ref{equ:3.5N}). Let $y_1,y_2\in Y\cap X_{r_i}$.
By the assumption of this theorem $G_{u_0}$ is transitive on $X_{r_i}$. 
Hence there exists $\sigma\in G_{u_0}$ satisfying
$\sigma(y_1)=y_2$. 
Since
$\sigma(u_0)=u_0$ and $u_s\in X_{l_s}$,
we must have $\sigma(u_s)\in X_{l_s}$ for $s=1,2,\ldots,N$.
Let
$\phi^{\sigma}_s=\phi^{(j_s)}_{\sigma(u_s)}$.
Since $\sigma(X_r)=X_r$ for $r=0,1,\ldots,d$ and
\begin{eqnarray}&&
\sum_{s=1}^Nc_s\phi^{\sigma}_s(x)
= \sum_{s=1}^dc_s\phi^{(j_s)}_{\sigma(u_s)}(x)
=\sum_{s=1}^Nc_s\frac{1}{|X|}E_{j_s}(x,\sigma(u_s))
\nonumber\\
&&
=\sum_{s=1}^Nc_s\frac{1}{|X|}E_{j_s}(\sigma^{-1}(x),u_s)
=\sum_{s=1}^Nc_s\phi_s(\sigma^{-1}(x)),
\end{eqnarray} 
$\{\phi^\sigma_1,\phi^\sigma_2,\ldots,\phi^\sigma_N\}$ 
is also a basis of
$L_0(S)+ L_1(S)+\cdots+ L_e(S)$.
Then we have
\begin{eqnarray}
&&
\langle \phi^\sigma_{ l_1},\phi^\sigma_{ l_2}\rangle
=\langle \phi^{(j_{l_1})}_{\sigma(u_{l_1})},
\phi^{(j_{l_2})}_{\sigma(u_{l_2})}\rangle
\nonumber\\
&&
=\sum_{i=1}^p
\frac{W_{r_i}}{|X_{r_i}|}
\sum_{\nu_1=0}^d\sum_{\nu_2=0}^d
\sum_{x\in X_{r_i}\cap\Gamma_{\nu_1}(\sigma(u_{l_1}))
\cap\Gamma_{\nu_2}(\sigma(u_{l_2}))}
\phi^{(j_{l_1})}_{\sigma(u_{l_1})}(x)
\phi^{(j_{l_2})}_{\sigma(u_{l_2})}(x)
\nonumber\\
&&
=\sum_{i=1}^p
\frac{W_{r_i}}{|X_{r_i}|}
\sum_{\nu_1=0}^d\sum_{\nu_2=0}^d
| \sigma(X_{r_i})\cap\Gamma_{\nu_1}(\sigma(u_{l_1}))
\cap\Gamma_{\nu_2}(\sigma(u_{l_2}))|
Q_{j_{l_1}}(\nu_1)Q_{j_{l_2}}(\nu_2)
\nonumber\\
&&
=\sum_{i=1}^p
\frac{W_{r_i}}{|X_{r_i}|}
\sum_{\nu_1=0}^d\sum_{\nu_2=0}^d
| X_{r_i}\cap\Gamma_{\nu_1}(u_{l_1})
\cap\Gamma_{\nu_2}(u_{l_2})|
Q_{j_{l_1}}(\nu_1)Q_{j_{l_2}}(\nu_2)
\nonumber\\
&&=\langle \phi^{(j_{l_1})}_{u_{l_1}},
\phi^{(j_{l_2})}_{u_{l_2}}\rangle
=\langle \phi_{l_1},\phi_{l_2}\rangle
\end{eqnarray}
for any $l_1,l_2\in \{1,\ldots,N\}$.
Here $\Gamma_\nu(u)=\{x\in X\mid (u,x)\in R_\nu\}$.
On the other hand if $y_1\in \Gamma_\nu(u_{l_1})$, 
then we must have 
$y_2=\sigma(y_1)\in \Gamma_\nu(\sigma(u_{l_1}))$.
Hence we have
\begin{eqnarray}\phi^\sigma_{l_1}(y_2)
=\phi^{(j_{l_1})}_{\sigma(u_{l_1})}(\sigma(y_1))
=\frac{1}{|X|}Q_{j_{l_1}}(\nu)
=\phi^{(j_{l_1})}_{u_{l_1}}(y_1)
=\phi_{l_1}(y_1).
\end{eqnarray}
Let $\{\varphi^{\sigma}_{1},\ldots,\varphi^{\sigma}_{N}\}$
be the orthonormal system obtained from 
$\{\phi^\sigma_1,\ldots,\phi^\sigma_N\}$
by the formulas (\ref{equ:3.5N}) and (\ref{equ:3.6N}).
Then we must have 
$\varphi_s(y_1)=\varphi^{\sigma}_{s}(y_2)$
for $s=1,2,\ldots, N$.
Hence we have
$$\sum_{s=1}^N(\varphi_s(y_1))^2
=\sum_{s=1}^N(\varphi^{\sigma}_s(y_2))^2.$$
This implies $w(y_1)=w(y_2)$ and completes the
proof of Theorem 1.1.
\hfill\qed\\


\section{Proof of Theorem \ref{theo:2.2}}
\subsection{Important propositions} 
It is known that the binary Hamming scheme 
$H(n,2)$ satisfies the assumption
of Theorem \ref{theo:2.1}.
First we introduce notation for $H(n,2)$.
Let $F=\{0,1\}$ and $X=F^n$
and $H(n,2)=(X,\{R_i\}_{0\leq i\leq n})$.
For $x=(x_1,x_2,\ldots,x_n)\in X$, we define 
$\overline{x}\subset
\{1,2,\ldots,n\}$ by 
$\overline{x}=\{i\mid x_i=1,~1\leq i\leq n\}$.

Let $(Y,w)$ be a relative tight $2e$-design of $H(n,2)$
supported by $p$ shells, i.e., $S=X_{r_1}\cup\cdots\cup X_{r_p}$.
Let $N=\dim(L_0(S)+L_1(S)+\cdots+L_e(S))$.
Then by Theorem \ref{theo:2.1} we have
$|Y|=N$ and $w(y)=w_{r_i}$ for any 
$y\in Y_{r_i}=Y\cap X_{r_i}$, $i=1,\ldots,p$, with
positive real numbers
$w_{r_1},\ldots,w_{r_p}$.
In the proof of Theorem \ref{theo:2.1}, we showed that
for any orthonormal basis $\{\varphi_1,\ldots,\varphi_N\}$
of $L_0(S)+L_1(S)+\cdots+L_e(S)$ with respect to the 
inner product defined by (\ref{equ:3.1}) then 
$$\sum_{i=1}^N\varphi_i(x)\varphi_i(y)
=\delta_{x,y}\frac{1}{w(y)}$$
holds holds for any $x,y\in Y$. 
We use this property and investigate the relations 
between the constants $N(=|Y|)$, 
$r_1,\ldots,r_p$, $N_{r_i}=|Y\cap X_{r_i}|$
$(1\leq i\leq p))$ and $w_{r_1},\ldots,w_{r_p}$. 
It is known that the first and second eigen matrices
$P$ and $Q$
of $H(n,2)$ coincide and given by
\begin{eqnarray}P_k(u)=Q_k(u)
=\sum_{i=0}^k(-1)^i{n-u\choose k-i}{u\choose i}.
\end{eqnarray}
In particular $k_i=m_i={n\choose i}$ holds for 
$i=0,1,\ldots,n$.
We consider the relative $2$-design $(Y,w)$ with respect to
$u_0$, on 
$S=X_{r_1}\cup X_{r_2}$.
Without loss of generality 
we may assume $u_0=(0,0,\ldots,0)$. 
Then $x\in X_r$ if and only if $|\overline{x}|=r$.
Let $X_1=\{u_1,\ldots,u_{n}\}$ (note that $k_1=m_1=n$
in this case). 
We use the following notation. 
\begin{eqnarray}&&
\phi_0(x)=\phi^{(0)}_{u_0}(x)=|X|E_0(x,u_0),\\
&&
\phi_j(x)=\phi^{(1)}_{u_j}(x)=|X|E_1(x,u_j)
\end{eqnarray}
for any $x\in X$.
By
Proposition 2.2 (2) (b) in \cite{Li-B-B-2013},
$\{\phi_0|_S,\phi_1|_S,\ldots, \phi_n|_S\}$ is an basis of
$L_0(S)+L_1(S)$, $S=X_{r_1}\cup X_{r_2}$,
for any integers $r_1,r_2$ satisfying $1\leq r_1< r_2\leq n-1$.
(The condition $(k,l)\not= (1,n-1)$ of Proposition 2.2 (2) (b) in 
\cite{Li-B-B-2013} is not correct. It should be $(k,l)\not= (0,n), (n,0)$.
So in our case we assume $1\leq r_1<r_2\leq n-1$, and then
$\dim(L_0(S)+L_1(S))=n+1$ holds. 
)
In this case the inner product $\langle f,g\rangle$, $f,g\in L_0(S)+L_1(S)$,
$S=X_{r_1}\cup X_{r_2}$ is defined by
$$\langle f,g\rangle=\frac{W_{r_1}}{|X_{r_1}|}
\sum_{x\in X_{r_1}}f(x)g(x)
+\frac{W_{r_2}}{|X_{r_2}|}
\sum_{x\in X_{r_2}}f(x)g(x).$$
By definition, $W_{r_i}=N_{r_i}w_{r_i}$ holds for $i=1,2$.
The following propositions play the important
role for the proof of Theorem \ref{theo:2.2}.
\begin{pro}$ $
\begin{enumerate}
\item
$\langle \phi_i,\phi_0\rangle
=\frac{(n-2)}{n}\left((n-2r_1)W_{r_1}+(n-2r_2)
W_{r_2}\right)$
for $1\leq i\leq n$.
\item
$\langle \phi_i,\phi_i\rangle=\sum_{\nu=1}^2\frac{W_{r_\nu}}{n}
\bigg(4(n-4)r_\nu^2-4n(n-4)r_\nu+n(n-2)^2\bigg)$
for $1\leq i\leq n$.
\item
$\langle \phi_i,\phi_j\rangle=
\langle \phi_1,\phi_2\rangle$

$=
\sum_{\nu=1}^2\frac{W_{r_\nu}}{n(n-1)}
\bigg(4(n^2-5n+8)r_\nu^2-4n(n^2-5n+8)r_\nu+n(n-1)(n-2)^2\bigg)
$
for any $1\leq i\not=j\leq n$.
\end{enumerate}
\label{pro:4.1}
\end{pro}
We use the following notation: $d_0=\langle \phi_1,\phi_0\rangle$,
$c_0=\langle \phi_1,\phi_1\rangle$
and $c_2=\langle \phi_1,\phi_2\rangle$.
\begin{pro} Let $h_1,h_2,\ldots, h_{n+1}$
be the orthogonal basis of 
$L_0(S)+L_1(S)$ obtained from the bases
$\{\phi_1,\phi_2,\ldots,\phi_n,\phi_0\}$ by 
Gram-Schmidt's method with this ordering.
Then we have the following formulas.
\begin{eqnarray}
&&h_1=\phi_1,\label{equ:4.4}\\
&&h_i=\phi_i-\frac{c_2}{c_0+(i-2)c_2}
\sum_{j=1}^{i-1}\phi_j
\qquad\mbox{for $i=2,3,\ldots,n$},
\label{equ:4.5}
\\
&&h_{n+1}=\phi_0-\frac{d_0}{c_0+(n-1)c_2}
\sum_{j=1}^n\phi_j,\label{equ:4.6}
\\
&&\|h_1\|^2=c_0,\label{equ:4.7}
\\
&&\|h_i\|^2
=\frac{(c_0-c_2)(c_0+(i-1)c_2)
}{(c_0+(i-2)c_2)},
\quad\mbox{for $i=2,\ldots,n$},\label{equ:4.8}\\
&&\|h_{n+1}\|^2
=W_{r_1}+W_{r_2}- \frac{nd_0^2}{c_0+(n-1)c_2}.
\label{equ:4.9}
\end{eqnarray}
\label{pro:4.2}
\end{pro}
\begin{pro} 
\begin{enumerate}
\item $2\leq N_{r_1},N_{r_2}\leq n-1$ holds and
\begin{eqnarray}&&
\frac{w_{r_2}}{w_{r_1}}=\frac{N_{r_1}r_1(n-N_{r_1})(n-r_1)}{r_2(N_{r_1}-1)(n+1-N_{r_1})(n-r_2)}.
\end{eqnarray}
\item If there exists an nonzero even integer $\alpha_\nu$ satisfying 
$2\leq \alpha_\nu \leq 2r_\nu$, and $x,y\in X_{r_\nu}~(\nu=1,2),$ with
$(x,y)\in R_{\alpha_\nu}$, then the following holds.
\begin{eqnarray}
&&\alpha_1=\frac{ 2(n-r_1)r_1N_{r_1}}
{n(N_{r_1}-1)},\\
&&\alpha_2=\frac{2(n-r_2)(n+1-N_{r_1})r_2}
{n(n-N_{r_1})}.
\end{eqnarray}

\item If there exists an even integer
$\gamma$ satisfying 
$(x,y)\in R_{\gamma}$, for $x\in X_{r_1}$ and $y\in X_{r_2}$, then the following holds. 
\begin{eqnarray}&&
\gamma=\frac{n(r_1+r_2)-2r_1r_2}{n}.
\end{eqnarray} 
\end{enumerate}
\label{pro:4.3}
\end{pro}
\subsection{Proof of the propositions}
\noindent
{\bf Proof of Proposition \ref{pro:4.1}}\\
(1):
\begin{eqnarray}&&
\langle \phi_{i},\phi_{0}\rangle
=\frac{W_{r_1}}{|X_{r_1}|}\sum_{x\in X_{r_1}}
\phi^{(1)}_{u_i}(x)\phi^{(0)}_{u_0}(x)
+\frac{W_{r_2}}{|X_{r_2}|}\sum_{x\in X_{r_2}}
\phi^{(1)}_{u_i}(x)\phi^{(0)}_{u_0}(x)
\nonumber\\
&&=\frac{W_{r_1}}{|X_{r_1}|}\sum_{\nu=0}^n
\sum_{x\in X_{r_1}\cap\Gamma_\nu(u_i)}
Q_1(\nu)
+\frac{W_{r_2}}{|X_{r_2}|}\sum_{x\in X_{r_2}}
\phi_{u_i}(x)
\nonumber\\
&&
=\frac{W_{r_1}}{|X_{r_1}|}\bigg(
|X_{r_1}\cap\Gamma_{r_1-1}(u_i)|
Q_1(r_1-1)
+|X_{r_1}\cap\Gamma_{r_1+1}(u_i)|
Q_1(r_1+1)
\bigg)\nonumber\\
&&
+\frac{W_{r_2}}{|X_{r_2}|}\bigg(
|X_{r_2}\cap\Gamma_{r_2-1}(u_i)|
Q_1(r_2-1)
+|X_{r_2}\cap\Gamma_{r_2+1}(u_i)|
Q_1(r_2+1)
\bigg)\nonumber
\\
&&=\sum_{\nu=1}^2\frac{W_{r_\nu}}{{n\choose r_\nu}}
\bigg({n-1\choose r_\nu-1}
Q_1(r_\nu-1)+{n-1\choose r_\nu}Q_1(r_\nu+1)
\bigg)
\nonumber\\
&&=\frac{(n-2)(n-2r_1)W_{r_1}}{n}+\frac{(n-2)(n-2r_2)W_{r_2}}{n}.
\end{eqnarray}
This proves Proposition \ref{pro:4.1} (1).\\
(2) and (3):\\
Let $1\leq i, j\leq n$. Then
\begin{eqnarray}&&
\langle \phi_{i},\phi_{j}\rangle
=\frac{W_{r_1}}{|X_{r_1}|}\sum_{x\in X_{r_1}}
\phi^{(1)}_{u_i}(x)\phi^{(1)}_{u_j}(x)
+\frac{W_{r_2}}{|X_{r_2}|}\sum_{x\in X_{r_2}}
\phi^{(1)}_{u_i}(x)\phi^{(1)}_{u_j}(x)
\nonumber
\\
&&=\frac{W_{r_1}}{|X_{r_1}|}
\sum_{l_1=0}^n\sum_{l_2=0}^n|X_{r_1}\cap
\Gamma_{l_1}(u_i)\cap\Gamma_{l_2}(u_j)|
Q_1(l_1)Q_1(l_2)
\nonumber
\\
&&
+\frac{W_{r_2}}{|X_{r_2}|}
\sum_{l_1=0}^n\sum_{l_2=0}^n|X_{r_2}\cap
\Gamma_{l_1}(u_i)\cap\Gamma_{l_2}(u_j)|
Q_1(l_1)Q_1(l_2).
\end{eqnarray}
If $u=u_i=u_j\in X_1$, then
$$|X_{r_i}\cap\Gamma_{\nu}(u)|=p_{r_i,\nu}^1
=\left\{\begin{array}{cl}
{n-1\choose r_i-1}&\mbox{if $\nu=r_i-1$},\\
{n-1\choose r_i+1}&\mbox{if $\nu=r_i+1$},\\
0&\mbox{otherwise}.
\end{array}\right.$$
Therefore we have
\begin{eqnarray}&&
\langle \phi^{(1)}_{u},\phi^{(1)}_{u}\rangle
=\sum_{\nu=1}^2\frac{W_{r_\nu}}{|X_{r_\nu}|}
\left({n-1\choose r_\nu-1}
Q_1(r_\nu-1)^2+{n-1\choose r_\nu}
Q_1(r_\nu+1)^2\right)
\nonumber\\
&&=\sum_{\nu=1}^2\frac{W_{r_\nu}}{{n\choose r_\nu}}
\left({n-1\choose r_\nu-1}
(n-2r_\nu+2)^2
+{n-1\choose r_\nu}
(n-2r_\nu-2)^2\right)
\nonumber\\
&&=\sum_{\nu=1}^2\frac{r_\nu!(n-r_\nu)!
W_{r_\nu}}
{n!}
\bigg(\frac{(n-1)!}{(r_\nu-1)!(n-r_\nu)!}
(n-2r_\nu+2)^2
\nonumber\\
&&+\frac{(n-1)!}{(r_\nu)!(n- r_\nu-1)!}
(n-2r_\nu-2)^2\bigg)
\nonumber\\
&&=\sum_{\nu=1}^2W_{r_\nu}
\bigg(\frac{r_\nu}{n}
(n-2r_\nu+2)^2
+\frac{(n-r_\nu)}{n}
(n-2r_\nu-2)^2\bigg)\nonumber\\
&&
=\sum_{\nu=1}^2\frac{W_{r_\nu}}{n}
\bigg(4(n-4)r_\nu^2-4n(n-4)r_\nu+n(n-2)^2\bigg).
\end{eqnarray}
This implies (2).\\
If $u_i\not=u_j$, and $1\leq r_1<r_2\leq n-1$, then
\begin{eqnarray}&&
\langle \phi_i,\phi_j\rangle
=\sum_{\nu=1}^2
\frac{W_{r_\nu}}{|X_{r_\nu}|}\sum_{x\in X_{r_\nu}}
\phi^{(1)}_{u_i}(x)\phi^{(1)}_{u_j}(x)\nonumber\\
&&
=\sum_{\nu=1}^2
\frac{W_{r_\nu}}{|X_{r_\nu}|}
\sum_{l_1,l_2=r_\nu-1,r_\nu+1}|X_{r_\nu}
\cap\Gamma_{l_1}(u_i)\cap\Gamma_{l_2}(u_j)|
Q_1(l_1)Q_1(l_2)
\nonumber
\\
&&
=\sum_{\nu=1}^2\frac{W_{r_\nu}}{n(n-1)}
\bigg(4(n^2-5n+8)r_\nu^2-4n(n^2-5n+8)r_\nu+n(n-1)(n-2)^2\bigg).
\nonumber\\
&&
\end{eqnarray}

\hfill\qed\\

\noindent
{\bf Proof of Proposition \ref{pro:4.2}}\\
(\ref{equ:4.4}) and (\ref{equ:4.7}) are already shown.
According to the Gram-Schmidt's method
let $h_i=\phi_i+\sum_{j=1}^{i-1}a_{i,j}\phi_j$
for $i=2,\ldots,n$.
Then $\langle h_1,h_2\rangle=0$ implies
$$h_2=\phi_2-\frac{c_2}{c_0}\phi_1.$$
Then we must have
\begin{eqnarray}&&\langle h_2,h_2\rangle
=\langle \phi_2,\phi_2\rangle
-2\frac{c_2}{c_0}\langle \phi_1,\phi_2\rangle
+(\frac{c_2}{c_0})^2
\langle\phi_1,\phi_1\rangle\nonumber\\
&&=c_0-2\frac{c_2^2}{c_0}+(\frac{c_2}{c_0})^2c_0
=\frac{(c_0-c_2)(c_0+c_2)}{c_0}.
\nonumber
\end{eqnarray}
Thus $h_2$ satisfies (\ref{equ:4.5}) and (\ref{equ:4.8}).
We prove (\ref{equ:4.5}) and (\ref{equ:4.8}) by induction on $i$.
Assume that (\ref{equ:4.5}) and (\ref{equ:4.8})
hold for any $i\leq s-1, s\leq n$ and we will show
that they also hold for $i=s$.
\begin{eqnarray}
&&
0=\langle h_s,h_{s-1}\rangle
=\langle\phi_s
+\sum_{j=1}^{s-1}a_{s,j}\phi_j
,h_{s-1}\rangle
=\langle\phi_s
+a_{s,s-1}\phi_{s-1}
,h_{s-1}\rangle\nonumber\\
&&=\langle\phi_s
+a_{s,s-1}\phi_{s-1}
,\phi_{s-1}\rangle
-\frac{c_2}{c_0+(s-3)c_2}\sum_{j=1}^{s-2}\langle\phi_s
+a_{s,s-1}\phi_{s-1}
,\phi_j\rangle\nonumber\\
&&=c_2+a_{s,s-1}c_0-\frac{c_2}{c_0+(s-3)c_2}
(s-2)(1+a_{s,s-1})c_2.
 \end{eqnarray}
This implies $a_{s,s-1}=-\frac{c_2}{c_0+(s-2)c_2}$.
By continuing such straight forwarded computation 
we obtain $a_{s,1}=a_{s,2}=\cdots=a_{s,s-1}
=-\frac{c_2}{c_0+(s-2)c_2}$ and we can verify the
formula (\ref{equ:4.8}) for $\|h_s\|^2$.
This completes the proof for 
(\ref{equ:4.5}) and  (\ref{equ:4.8}).
Next let
$$h_{n+1}=\phi_0+\sum_{j=1}^na_j\phi_j.$$
Then we have the following.
\begin{eqnarray}&&
0=\langle h_{n+1}, h_n\rangle=
\langle \phi_0+\sum_{j=1}^na_j\phi_j, h_n\rangle=
\langle \phi_0+a_n\phi_n, h_n\rangle
\nonumber\\
&&=\langle \phi_0+a_n\phi_n, \phi_n
-\frac{c_2}{c_0+(n-2)c_2}
\sum_{j=1}^{n-1}\phi_j\rangle
\nonumber\\
&&=\langle \phi_0, \phi_n\rangle
+a_n\langle \phi_n, \phi_n\rangle
-\frac{c_2}{c_0+(n-2)c_2}\langle \phi_0, 
\sum_{j=1}^{n-1}\phi_j\rangle
\nonumber\\
&&-a_n\frac{c_2}{c_0+(n-2)c_2}\langle \phi_n, 
\sum_{j=1}^{n-1}\phi_j\rangle
\nonumber\\
&&=d_0(1
-\frac{(n-1)c_2}{c_0+(n-2)c_2})
+(c_0-\frac{(n-1)c_2^2}{c_0+(n-2)c_2})a_n.
\end{eqnarray}
Hence 
\begin{eqnarray}&&a_n=-\frac{d_0}{c_0+(n-1)c_2}
\nonumber\\
&&=
-\frac{(n-2r_1)W_{r_1}+(n-2r_2)W_{r_2}}
{(n-2)((n-2r_1)^2W_{r_1}+(n-2r_2)^2W_{r_2})}.
\end{eqnarray}
By straight forwarded computation we obtain
$a_1=a_2=\cdots=a_n$.   
This implies (\ref{equ:4.6}).
\begin{eqnarray}&&
\|h_{n+1}\|^2=\langle h_{n+1},h_{n+1}\rangle
\nonumber\\
&&
=
\langle \phi_0-\frac{d_0}{c_0+(n-1)c_2}
\sum_{j=1}^{n}\phi_j,
\phi_0-\frac{d_0}{c_0+(n-1)c_2}
\sum_{l=1}^{n}\phi_l\rangle
\nonumber\\
&&=\langle \phi_0,
\phi_0\rangle
-\frac{2d_0}{c_0+(n-1)c_2}
\sum_{l=1}^{n}\langle \phi_0,\phi_l\rangle
+ \left(\frac{d_0}{c_0+(n-1)c_2}\right)^2\sum_{j=1}^{n}
\langle\phi_j,\sum_{l=1}^{n}\phi_l\rangle
\nonumber
\end{eqnarray}
\begin{eqnarray}&&=W_{r_1}+W_{r_2}-\frac{2nd_0^2}{c_0+(n-1)c_2}
+ \left(\frac{d_0}{c_0+(n-1)c_2}\right)^2n(c_0+(n-1)c_2)\nonumber\\
&&=W_{r_1}+W_{r_2}- \frac{nd_0^2}{c_0+(n-1)c_2}.
\end{eqnarray}
This completes the proof of Proposition \ref{pro:4.2}.
\hfill\qed\\

\noindent
{\bf Proof of Proposition \ref{pro:4.3}}\\
(1): Let $\nu=1$ or $2$.
Let $x\in Y_{r_\nu}=X_{r_\nu}\cap Y$. Choose the ordering of the elements 
in $X_1$, we may assume $\overline{x}=\{1,2,\ldots,r_\nu\}$.
Then 
$\phi_i(x)=Q_1(r_\nu-1)$ for $i=1,2,\ldots,r_\nu$
and $\phi_i(x)=Q_1(r_\nu+1)$ for $i=r_\nu+1,r_\nu+2,\ldots,n$.
Hence 
\begin{eqnarray}&&h_i(x)=Q_1(r_\nu-1)-\frac{(i-1)c_2}{c_0+(i-2)c_2}Q_1(r_\nu-1) 
~\mbox{for $1\leq i\leq r_\nu $},\\
&&h_i(x)=Q_1(r_\nu+1)-\frac{r_\nu c_2}{c_0+(i-2)c_2}Q_1(r_\nu-1) 
-\frac{(i-r_\nu-1)c_2}{c_0+(i-2)c_2}Q_1(r_\nu+1) 
\nonumber\\
&&
\mbox{for $r_\nu+1\leq i\leq n $},\\
&&h_{n+1}(x)=1-\frac{d_0}{c_0+(n+1-2)c_2}
(r_\nu Q_1(r_\nu-1)+(n-r_\nu)Q_1(r_\nu+1)).
\nonumber\\
&&
\end{eqnarray}
Let $\varphi_i=\frac{1}{\|h_i\|}h_i$, $i=1,2,\ldots,n+1$.
Then $\{\varphi_1,\ldots,\varphi_{n+1}\}$ is an orthnormal 
basis of $L_0(S)+L_1(S)$. Hence we have
{\small
\begin{eqnarray}&&
\frac{1}{w_\nu}=\sum_{i=1}^{n+1}\varphi_i(x)^2=
\sum_{i=1}^{n+1}\frac{1}{\|h_i\|^2}h_i(x)^2
\nonumber\\
&&=\sum_{s=1}^{r_\nu}
\frac{c_0+(s-2)c_2}{(c_0-c_2)(c_0+(s-1)c_2)}
\bigg(1-\frac{(s-1)c_2}{(c_0+(s-2)c_2)}
\bigg)^2 Q_1(r_\nu-1)^2
\nonumber\\
&&
+\sum_{s=r_\nu+1}^n\frac{c_0+(s-2)c_2}{(c_0-c_2)(c_0+(s-1)c_2)}
\bigg(Q_1(r_\nu+1)
-\frac{r_\nu c_2}{c_0+(s-2)c_2}Q_1(r_\nu-1)
\nonumber\\
&&
-\frac{(s-r_\nu-1)c_2}{c_0+(s-2)c_2}Q_1(r_\nu+1)\bigg)^2
\nonumber\\
&&
+\frac{c_0+(n-1)c_2}{(W_{r_1}+W_{r_2})
(c_0+(n-1)c_2)
-nd_0^2}\bigg\{1-\frac{d_0}{c_0+(n-1)c_2} 
\bigg(r_\nu Q_1(r_\nu-1)
\nonumber\\
&&
+(n-r_\nu)Q_1(r_\nu+1)\bigg)\bigg\}^2.
\label{equ:4.25}
\end{eqnarray}
}
{\small
\begin{eqnarray}
&&\sum_{s=1}^{r_\nu}
\frac{c_0+(s-2)c_2}{(c_0-c_2)(c_0+(s-1)c_2)}
\bigg(1-\frac{(s-1)c_2}{(c_0+(s-2)c_2)}
\bigg)^2 Q_1(r_\nu-1)^2
\nonumber\\
&&=\frac{(c_0-c_2)(n-2r_\nu+2)^2}{c_2}
\sum_{s=1}^{r_\nu}\left(
\frac{1}{c_0+(s-2)c_2}-\frac{1}{c_0+(s-1)c_2}
\right)
\nonumber\\
&&=\frac{(c_0-c_2)(n-2r_\nu+2)^2}{c_2}
\left(
\frac{1}{c_0-c_2}-\frac{1}{c_0+(r_\nu-1)c_2}
\right).\label{equ:4.26}
\end{eqnarray}
}
{\small
\begin{eqnarray}
&&\sum_{s=r_\nu+1}^n\frac{c_0+(s-2)c_2}{(c_0-c_2)(c_0+(s-1)c_2)}
\bigg(Q_1(r_\nu+1)
-\frac{r_\nu c_2}{c_0+(s-2)c_2}Q_1(r_\nu-1)
\nonumber\\
&&
-\frac{(s-r_\nu-1)c_2}{c_0+(s-2)c_2}Q_1(r_\nu+1)\bigg)^2
\nonumber\\
&&=\frac{(nc_0-nc_2-2r_\nu c_0-2r_\nu c_2-2c_0+2c_2)^2}{(c_0-c_2)c_2}\sum_{s=r_\nu+1}^n
\left(\frac{1}{c_0+(s-2)c_2}-
\frac{1}{c_0+(s-1)c_2}\right)
\nonumber\\
&&=\frac{(nc_0-nc_2-2r_\nu c_0-2r_\nu c_2-2c_0+2c_2)^2}{(c_0-c_2)c_2}
\left(\frac{1}{c_0+(r_\nu-1)c_2}-
\frac{1}{c_0+(n-1)c_2}\right)
\nonumber\\
&&=\frac{(nc_0-nc_2-2r_\nu c_0-2r_\nu c_2-2c_0+2c_2)^2(n-r_\nu)}{(c_0-c_2)(c_0+r_\nu c_2-c_2)(c_0+nc_2-c_2))}.
\label{equ:4.27}
\end{eqnarray}
}
{\small
\begin{eqnarray}&&
\frac{c_0+(n-1)c_2}{(W_{r_1}+W_{r_2})
(c_0+(n-1)c_2)
-nd_0^2}\bigg\{1-\frac{d_0}{c_0+(n-1)c_2} 
\bigg(r_\nu Q_1(r_\nu-1)
\nonumber\\
&&
+(n-r_\nu)Q_1(r_\nu+1)\bigg)\bigg\}^2=
\frac{(c_0+nc_2-c_2+2d_0r_\nu n-4d_0r_\nu-d_0n^2+2d_0n)^2}{(c_0+(n-1)c_2)
((W_{r_1}+W_{r_2})(c_0+(n-1)c_2)-nd_0^2)}.
\nonumber\\
&&\label{equ:4.28}
\end{eqnarray}
}
Since $W_{r_1}=N_{r_1}w_{r_1}$ and
$W_{r_2}=(n+1-N_{r_1})w_{r_2}$,
(\ref{equ:4.25}), (\ref{equ:4.26}),  (\ref{equ:4.27}) and (\ref{equ:4.28}) imply 
\begin{eqnarray}&&
\frac{1}{w_{r_2}}
=\frac{nr_2(n+1-N_{r_1})(n-r_2)w_{r_2}+N_{r_1}r_1(n-r_1)w_{r_1}}
{(n+1-N_{r_1})\bigg(N_{r_1}w_{r_1}r_1(n-r_1)+r_2(n+1-N_{r_1})w_{r_2}(n-r_2)\bigg)w_{r_2}}.
\nonumber\\
&&
\end{eqnarray}
Therefore we have
$$w_{r_2}=\frac{N_{r_1}r_1(n-N_{r_1})(n-r_1)}{r_2(N_{r_1}-1)(n+1-N_{r_1})(n-r_2)}w_{r_1}.$$
This completes the proof for (1).\\
\noindent
(2): Let $\nu=1$ or $2$.
Let $x,y\in Y_{r_\nu}$ and $x\not= y$.
then $(x,y)\in R_{\alpha_\nu}$, $\alpha_\nu=2,\ldots,2r_\nu$.
Let $(x,y)\in R_{\alpha_\nu}$.
Then take the ordering of $\{u_1,\ldots,u_n\}$
so that $\overline{x}=\{1,2,\ldots,r_\nu\}$ and
$\overline{y}=\{\frac{1}{2}\alpha_\nu+1,
\frac{1}{2}\alpha_\nu+2,\cdots,\frac{1}{2}\alpha_\nu
+r_\nu\}$.
Then
for $1\leq i\leq \frac{1}{2}\alpha_\nu$ we have
\begin{eqnarray}%
&&h_i(x)=Q_1(r_\nu-1)-\frac{(i-1)c_2}{c_0+(i-2)c_2}Q_1(r_\nu-1),\\
&&h_i(y)=Q_1(r_\nu+1)-\frac{(i-1)c_2}{c_0+(i-2)c_2}Q_1(r_\nu+1),
\end{eqnarray}
If $\frac{1}{2}\alpha_\nu+1\leq i\leq r_\nu$, then
\begin{eqnarray}
&&h_i(x)=Q_1(r_\nu-1)-\frac{(i-1)c_2}{c_0+(i-2)c_2}Q_1(r_\nu-1),\\
&&h_i(y)=Q_1(r_\nu-1)-\frac{\frac{1}{2}\alpha_\nu c_2}{c_0+(i-2)c_2}Q_1(r_\nu+1)
-\frac{(i-\frac{1}{2}\alpha_\nu-1) c_2}{c_0+(i-2)c_2}Q_1(r_\nu-1).\nonumber\\
&&
\end{eqnarray}
If $r_\nu+1\leq i\leq r_\nu+\frac{1}{2}\alpha_\nu$, then
\begin{eqnarray}
&&h_i(x)=Q_1(r_\nu+1)
-\frac{r_\nu c_2}{c_0+(i-2)c_2}Q_1(r_\nu-1)
-\frac{(i-r_\nu-1)c_2}{c_0+(i-2)c_2}Q_1(r_\nu+1),
\nonumber\\
&&h_i(y)=Q_1(r_\nu-1)-\frac{\frac{1}{2}\alpha_\nu c_2}{c_0+(i-2)c_2}Q_1(r_\nu+1)
-\frac{(i-\frac{1}{2}\alpha_\nu-1) c_2}{c_0+(i-2)c_2}Q_1(r_\nu-1),\nonumber\\
&&
\end{eqnarray}
If $r_\nu+\frac{1}{2}\alpha_\nu+1\leq i\leq n$, then
\begin{eqnarray}
&&h_i(x)=Q_1(r_\nu+1)
-\frac{r_\nu c_2}{c_0+(i-2)c_2}Q_1(r_\nu-1)
-\frac{(i-r_\nu-1)c_2}{c_0+(i-2)c_2}Q_1(r_\nu+1),
\nonumber\\
&&\\
&&h_i(y)=Q_1(r_\nu+1)
-\frac{r_\nu c_2}{c_0+(i-2)c_2}Q_1(r_\nu-1)
-\frac{(i-r_\nu-1) c_2}{c_0+(i-2)c_2}Q_1(r_\nu+1).
\nonumber\\
&&
\end{eqnarray}
If $i=n+1$, then
\begin{eqnarray}
&&h_{n+1}(x)=h_{n+1}(y)
=1-\frac{d_0}{c_0+(n-1)c_2}\bigg(r_\nu Q_1(r_\nu-1)
+(n-r_\nu)Q_1(r_\nu+1)\bigg).\nonumber\\
&&
\end{eqnarray}
Then
{\small
\begin{eqnarray}
&&\sum_{i=1}^{n+1}\varphi_i(x)\varphi_i(y)
=\sum_{i=1}^{n+1}\frac{1}{\|h_i\|^2}h_i(x)h_i(y)\nonumber\\
&&
=\sum_{i=1}^{\frac{1}{2}\alpha_\nu}
\frac{c_0+(i-2)c_2}{(c_0-c_2)(c_0+(i-1)c_2)}
\bigg(Q_1(r_\nu-1)-\frac{(i-1)c_2}{c_0+(i-2)c_2}Q_1(r_\nu-1)
\bigg)\times\nonumber\\
&&\bigg(Q_1(r_\nu+1)-\frac{(i-1)c_2}{c_0+(i-2)c_2}
Q_1(r_\nu+1)
\bigg)\nonumber\\
&&
+\sum_{i=\frac{1}{2}\alpha_\nu+1}^{r_\nu}\frac{c_0+(i-2)c_2}
{(c_0-c_2)(c_0+(i-1)c_2)}
\bigg(Q_1(r_\nu-1)-\frac{(i-1)c_2}{c_0+(i-2)c_2}Q_1(r_\nu-1)\bigg)\times\nonumber\\
&&\bigg(Q_1(r_\nu-1)-\frac{\frac{1}{2}\alpha_\nu c_2}{c_0+(i-2)c_2}Q_1(r_\nu+1)
-\frac{(i-\frac{1}{2}\alpha_\nu-1) c_2}{c_0+(i-2)c_2}Q_1(r_\nu-1)\bigg)\nonumber
\end{eqnarray}
\begin{eqnarray}
&&+\sum_{i=r_\nu+1}^{r_\nu+\frac{1}{2}\alpha_\nu}\frac{c_0+(i-2)c_2}{(c_0-c_2)(c_0+(i-1)c_2)}\times\nonumber\\
&&\bigg(Q_1(r_\nu+1)
-\frac{r_\nu c_2}{c_0+(i-2)c_2}Q_1(r_\nu-1)
-\frac{(i-r_\nu-1)c_2}{c_0+(i-2)c_2}Q_1(r_\nu+1)\bigg)\times
\nonumber\\
&&\bigg(Q_1(r_\nu-1)-\frac{\frac{1}{2}\alpha_\nu c_2}{c_0+(i-2)c_2}Q_1(r_\nu+1)
-\frac{(i-\frac{1}{2}\alpha_\nu-1) c_2}{c_0+(i-2)c_2}Q_1(r_\nu-1)\bigg)\nonumber\\
&&+\sum_{i=r_\nu+\frac{1}{2}\alpha_\nu+1}^{n}
\frac{c_0+(i-2)c_2}{(c_0-c_2)(c_0+(i-1)c_2)}\times\nonumber\\
&&\bigg(Q_1(r_\nu+1)
-\frac{r_\nu c_2}{c_0+(i-2)c_2}Q_1(r_\nu-1)
-\frac{(i-r_\nu-1)c_2}{c_0+(i-2)c_2}Q_1(r_\nu+1)\bigg)^2
\nonumber\\
&&\nonumber\\
&&+\frac{c_0+(n-1)c_2}
{(W_{r_1}+W_{r_2})(c_0+(n-1)c_2)-nd_0^2}\times\nonumber\\
&&
\bigg\{1-\frac{d_0}{c_0+(n-1)c_2}\bigg(r_\nu Q_1(r_\nu-1)
+(n-r_\nu)Q_1(r_\nu+1)\bigg)\bigg\}^2.
\label{equ:4.38}
\end{eqnarray}}
\begin{eqnarray}
&&\sum_{i=1}^{\frac{1}{2}\alpha_\nu}
\frac{c_0+(i-2)c_2}{(c_0-c_2)(c_0+(i-1)c_2)}
\bigg(Q_1(r_\nu-1)-\frac{(i-1)c_2}{c_0+(i-2)c_2}Q_1(r_\nu-1)
\bigg)\times\nonumber\\
&&\bigg(Q_1(r_\nu+1)-\frac{(i-1)c_2}{c_0+(i-2)c_2}
Q_1(r_\nu+1)
\bigg)\nonumber\\
&&=(c_0-c_2)Q_1(r_\nu-1)Q_1(r_\nu+1)
\sum_{i=1}^{\frac{1}{2}\alpha_\nu}
\frac{1}{(c_0+(i-2)c_2)(c_0+(i-1)c_2)}
\nonumber\\
&&=\frac{(c_0-c_2)Q_1(r_\nu-1)Q_1(r_\nu+1)}
{c_2}\sum_{i=1}^{\frac{1}{2}\alpha_\nu}\bigg(
\frac{1}{c_0+(i-2)c_2}-\frac{1}{c_0+(i-1)c_2}
\bigg)
\nonumber\\
&&=\frac{(c_0-c_2)(n-2r_\nu+2)(n-2r_\nu-2)}
{c_2}\bigg(
\frac{1}{c_0-c_2}-\frac{1}{c_0+(\frac{\alpha_\nu}{2}-1)c_2}
\bigg).\label{equ:4.39}
\end{eqnarray}
\begin{eqnarray}
&&\sum_{i=\frac{1}{2}\alpha_\nu+1}^{r_\nu}\frac{c_0+(i-2)c_2}
{(c_0-c_2)(c_0+(i-1)c_2)}
\bigg(Q_1(r_\nu-1)-\frac{(i-1)c_2}{c_0+(i-2)c_2}Q_1(r_\nu-1)\bigg)\times\nonumber\\
&&\bigg(Q_1(r_\nu-1)-\frac{\frac{1}{2}\alpha_\nu c_2}{c_0+(i-2)c_2}Q_1(r_\nu+1)
-\frac{(i-\frac{1}{2}\alpha_\nu-1) c_2}{c_0+(i-2)c_2}Q_1(r_\nu-1)\bigg)\nonumber\\
&&=
(n-2(r_\nu-1))((c_0-c_2)(n-2(r_\nu-1))+2\alpha_\nu c_2)\times\nonumber\\
&&\sum_{i=\frac{1}{2}\alpha_\nu+1}^{r_\nu}
\frac{1}{(c_0+(i-2)c_2)(c_0+(i-1)c_2)}
\nonumber\\
&&=\frac{(n-2(r_\nu-1))((c_0-c_2)(n-2(r_\nu-1))+2\alpha_\nu c_2)}{c_2}
\left(\frac{1}{c_0+(\frac{\alpha_\nu}{2}-1)c_2}-\frac{1}{c_0+(r_\nu-1)c_2}\right)
\nonumber
\end{eqnarray}
\begin{eqnarray}
&&
=\frac{(n-2r_\nu+2)(c_0(n-2r_\nu+2)-c_2(n-2r_\nu+2)+2\alpha_\nu c_2)(2r_\nu-\alpha_\nu)}
{(2c_0+(\alpha_\nu-2)c_2)(c_0+(r_\nu -1)c_2)}.
\label{equ:4.40}
\end{eqnarray}
\begin{eqnarray}
&&\sum_{i=r_\nu+1}^{r_\nu+\frac{1}{2}\alpha_\nu}\frac{c_0+(i-2)c_2}{(c_0-c_2)(c_0+(i-1)c_2)}\times\nonumber\\
&&\bigg(Q_1(r_\nu+1)
-\frac{r_\nu c_2}{c_0+(i-2)c_2}Q_1(r_\nu-1)
-\frac{(i-r_\nu-1)c_2}{c_0+(i-2)c_2}Q_1(r_\nu+1)\bigg)\times
\nonumber\\
&&\bigg(Q_1(r_\nu-1)-\frac{\frac{1}{2}\alpha_\nu c_2}{c_0+(i-2)c_2}Q_1(r_\nu+1)
-\frac{(i-\frac{1}{2}\alpha_\nu-1) c_2}{c_0+(i-2)c_2}Q_1(r_\nu-1)\bigg)\nonumber\\
&&
=\frac{(c_0(n-2r_\nu+2)
-c_2(n-2r_\nu+2)+2\alpha_\nu c_2)(c_0(n-2r_\nu-2)-c_2(n+2r_\nu-2))}{c_2(c_0-c_2)}
\times\nonumber\\
&&
\sum_{i=r_\nu+1}^{r_\nu+\frac{1}{2}\alpha_\nu}\left(\frac{1}{c_0+(i-2)c_2)}-\frac{1}{c_0+(i-1)c_2}\right)
\nonumber\\
&&
=\frac{(c_0(n-2r_\nu+2)
-c_2(n-2r_\nu+2)+2\alpha_\nu c_2)(c_0(n-2r_\nu-2)-c_2(n+2r_\nu-2))}{c_2(c_0-c_2)}
\times \nonumber\\
&&\left(\frac{1}{c_0+(r_\nu-1)c_2)}
-\frac{1}{c_0+(r_\nu+\frac{\alpha_\nu}{2}-1)c_2}\right).
\label{equ:4.41}
\end{eqnarray}
\begin{eqnarray}
&&\sum_{i=r_\nu+\frac{1}{2}\alpha_\nu+1}^{n}
\frac{c_0+(i-2)c_2}{(c_0-c_2)(c_0+(i-1)c_2)}\times\nonumber\\
&&\bigg(Q_1(r_\nu+1)
-\frac{r_\nu c_2}{c_0+(i-2)c_2}Q_1(r_\nu-1)
-\frac{(i-r_\nu-1)c_2}{c_0+(i-2)c_2}Q_1(r_\nu+1)\bigg)^2
\nonumber\\
&&=\frac{(c_0(n-2r_\nu-2)-c_2(n+2r_\nu-2))^2}{(c_0-c_2)c_2}
\left(\frac{1}{c_0+(i-2)c_2}-\frac{1}{c_0
+(i-1)c_2}\right)
\nonumber\\
&&=\frac{(c_0(n-2r_\nu-2)-c_2(n+2r_\nu-2))^2}{(c_0-c_2)c_2}
\left(\frac{1}{c_0+(r_\nu+\frac{1}{2}\alpha_\nu-1)c_2}-\frac{1}{c_0
+(n-1)c_2}\right).
\nonumber\\
&&\label{equ:4.42}
\end{eqnarray}
(\ref{equ:4.38}),
(\ref{equ:4.39}), (\ref{equ:4.40}), (\ref{equ:4.41}), and (\ref{equ:4.42}) imply
\begin{eqnarray}
&&-\frac{8\alpha_\nu}{c_0-c_2}
+\frac{c_0(4(n-4)r_\nu^2-4n(n-4)r_\nu
+n(n-2)^2)-nc_2(4r_\nu^2-4n r_\nu+(n-2)^2)}{
(c_0+(n-1)c_2)(c_0-c_2)}
\nonumber\\
&&
+\frac{c_0+(n-1)c_2}
{(W_{r_1}+W_{r_2})(c_0+(n-1)c_2)-nd_0^2}
\bigg\{1-\frac{(n-2)(n-2r_\nu)d_0}{c_0+(n-1)c_2}
\bigg\}^2=0.
\label{equ:4.43}
\end{eqnarray}
Then using the formula in Proposition \ref{pro:4.1}, we have
$$\frac{2nW_{r_1}r_1(n-r_1)+2W_{r_2}r_2(n-r_2)-n\alpha_\nu W_{r_\nu}(n-1)}{2W_{r_\nu}(W_{r_2}r_2(n-r_2)+W_{r_1}r_1(n-r_1))}=0.
$$
This implies
$$\alpha_1=\frac{2(
n(n-r_1)r_1W_{r_1}
+(n-r_2)r_2W_{r_2})}{n(n-1)W_{r_1}}
$$
and
$$\alpha_2=\frac{2(
n(n-r_2)r_2W_{r_2}
+(n-r_1)r_1W_{r_1})}{n(n-1)W_{r_2}}.
$$
Then substitute $W_{r_1}=N_{r_1}w_{r_1}$,
$W_{r_2}=(n+1-N_{r_1})w_{r_2}$,
where $w_{r_2}$ is given in (1), 
we obtain (2).\\
\noindent
(3): Let $x\in Y_{r_1}$ and $y\in Y_{r_2}$.
Then $(x,y)\in R_{r_2-r_1+2a}$ with an integer
$a$ satisfying $0\leq a\leq r_1$.
Choose the ordering of the elements in $X_1$ so that
$\overline{x}=\{1,2,\ldots,r_1\}$ and 
$\overline{y}=\{a+1,
a+2,\ldots,a+r_2\}$.
hold. 

\noindent
Then for $1\leq i\leq a$ we have
\begin{eqnarray}%
&&h_i(x)=Q_1(r_1-1)-\frac{(i-1)c_2}{c_0+(i-2)c_2}Q_1(r_1-1),\\
&&h_i(y)=Q_1(r_2+1)-\frac{(i-1)c_2}{c_0+(i-2)c_2}Q_1(r_2+1).
\end{eqnarray}
If $a+1\leq i\leq r_1$, then
\begin{eqnarray}
&&h_i(x)=Q_1(r_1-1)-\frac{(i-1)c_2}{c_0+(i-2)c_2}Q_1(r_1-1),\\
&&h_i(y)=Q_1(r_2-1)-\frac{a c_2}{c_0+(i-2)c_2}Q_1(r_2+1)
-\frac{(i-a-1) c_2}{c_0+(i-2)c_2}Q_1(r_2-1).\nonumber\\
&&
\end{eqnarray}
If $r_1+1\leq i\leq r_2+a$, then
\begin{eqnarray}
&&h_i(x)=Q_1(r_1+1)-\frac{r_1c_2}{c_0+(i-2)c_2}Q_1(r_1-1)
-\frac{(i-r_1-1)c_2}{c_0+(i-2)c_2}Q_1(r_1+1),\nonumber\\
&&\\
&&h_i(y)=Q_1(r_2-1)-\frac{a c_2}{c_0+(i-2)c_2}Q_1(r_2+1)
-\frac{(i-a-1)c_2}{c_0+(i-2)c_2}Q_1(r_2-1).\nonumber\\
&&
\end{eqnarray}
If $r_2+a\leq i\leq n$, then
\begin{eqnarray}
&&h_i(x)=Q_1(r_1+1)-\frac{r_1c_2}{c_0+(i-2)c_2}Q_1(r_1-1)
-\frac{(i-r_1-1)c_2}{c_0+(i-2)c_2}Q_1(r_1+1),\nonumber\\
&&\\
&&h_i(y)=Q_1(r_2+1)-\frac{r_2c_2}{c_0+(i-2)c_2}Q_1(r_2-1)
-\frac{(i-r_2-1)c_2}{c_0+(i-2)c_2}Q_1(r_2+1).\nonumber\\
\end{eqnarray}
If $ i=n+1$, then
\begin{eqnarray}
&&h_{n+1}(x)
=1-\frac{d_0}{c_0+(n-1)c_2}\bigg(r_1Q_1(r_1-1)
+(n-r_1)Q_1(r_1+1)\bigg),\nonumber\\
&&\\
&&h_{n+1}(y)=1-\frac{d_0}{c_0+(n-1)c_2}\bigg(r_2Q_1(r_2-1)
+(n-r_2)Q_1(r_2+1)\bigg).\nonumber\\
&&
\end{eqnarray}
Hence we have
{\small
\begin{eqnarray}
&&\sum_{i=1}^{n+1}\varphi_i(x)\varphi_i(y)\nonumber\\
&&=\sum_{i=1}^{a}\frac{c_0+(i-2)c_2}{(c_0-c_2)(c_0+(i-1)c_2)}
\bigg(Q_1(r_1-1)-\frac{(i-1)c_2}{c_0+(i-2)c_2}Q_1(r_1-1)\bigg)\times
\nonumber\\
&&
\bigg(Q_1(r_2+1)-\frac{(i-1)c_2}{c_0+(i-2)c_2}Q_1(r_2+1)\bigg)
\nonumber\\
&&+\sum_{i=a+1}^{r_1}\frac{c_0+(i-2)c_2}{(c_0-c_2)(c_0+(i-1)c_2)}
\bigg(Q_1(r_1-1)-\frac{(i-1)c_2}{c_0+(i-2)c_2}Q_1(r_1-1)\bigg)\times
\nonumber\\
&&
\bigg(Q_1(r_2-1)-\frac{a c_2}{c_0+(i-2)c_2}Q_1(r_2+1)
-\frac{(i-a-1) c_2}{c_0+(i-2)c_2}Q_1(r_2-1)\bigg)\nonumber\\
&&+\sum_{i=r_1+1}^{r_2+a}\frac{c_0+(i-2)c_2}{(c_0-c_2)(c_0+(i-1)c_2)}
\times
\nonumber\\
&&
\bigg(Q_1(r_1+1)-\frac{r_1c_2}{c_0+(i-2)c_2}Q_1(r_1-1)
-\frac{(i-r_1-1)c_2}{c_0+(i-2)c_2}Q_1(r_1+1)\bigg)\times
\nonumber\\
&&\bigg(Q_1(r_2-1)-\frac{a c_2}{c_0+(i-2)c_2}Q_1(r_2+1)
-\frac{(i-a-1)c_2}{c_0+(i-2)c_2}Q_1(r_2-1)\bigg)\nonumber\\
&&+\sum_{i=r_2+a+1}^{n}\frac{c_0+(i-2)c_2}{(c_0-c_2)(c_0+(i-1)c_2)}
\times
\nonumber\\
&&\bigg(Q_1(r_1+1)-\frac{r_1c_2}{c_0+(i-2)c_2}Q_1(r_1-1)
-\frac{(i-r_1-1)c_2}{c_0+(i-2)c_2}Q_1(r_1+1)\bigg)\times
\nonumber\\
&&\bigg(Q_1(r_2+1)-\frac{r_2c_2}{c_0+(i-2)c_2}Q_1(r_2-1)
-\frac{(i-r_2-1)c_2}{c_0+(i-2)c_2}Q_1(r_2+1)\bigg)\nonumber\\
&&+\frac{c_0+(n-1)c_2}
{(W_{r_1}+W_{r_2})(c_0+(n-1)c_2)-nd_0^2}\times\nonumber\\
&&\bigg\{
1-\frac{d_0}{c_0+(n-1)c_2}\bigg(r_1Q_1(r_1-1)
+(n-r_1)Q_1(r_1+1)\bigg)\bigg\}\times\nonumber\\
&&\bigg\{1-\frac{d_0}{c_0+(n-1)c_2}\bigg(r_2Q_1(r_2-1)
+(n-r_2)Q_1(r_2+1)\bigg)\bigg\}.
\label{equ:N4.54}
\end{eqnarray}
}

\begin{eqnarray}
&&\sum_{i=1}^{a}\frac{c_0+(i-2)c_2}{(c_0-c_2)(c_0+(i-1)c_2)}
\bigg(Q_1(r_1-1)-\frac{(i-1)c_2}{c_0+(i-2)c_2}Q_1(r_1-1)\bigg)\times
\nonumber\\
&&
\bigg(Q_1(r_2+1)-\frac{(i-1)c_2}{c_0+(i-2)c_2}Q_1(r_2+1)\bigg)
\nonumber\\
&&=(c_0-c_2)(n-2(r_1-1))(n-2(r_2+1))\sum_{i=1}^{a}
\frac{1}{(c_0+(i-2)c_2)(c_0+(i-1)c_2)}
\nonumber\\
&&=\frac{1}{c_2}(c_0-c_2)(n-2(r_1-1))(n-2(r_2+1))\left(
\frac{1}{c_0-c_2}-\frac{1}{c_0+(a-1)c_2}\right)
\nonumber\\
&&
=\frac{(n-2(r_1-1))(n-2(r_2+1))a}{(c_0-c_2)+ac_2}.
\label{equ:N4.55}
\end{eqnarray}
\begin{eqnarray}
&&\sum_{i=a+1}^{r_1}\frac{c_0+(i-2)c_2}{(c_0-c_2)(c_0+(i-1)c_2)}
\bigg(Q_1(r_1-1)-\frac{(i-1)c_2}{c_0+(i-2)c_2}Q_1(r_1-1)\bigg)\times
\nonumber\\
&&
\bigg(Q_1(r_2-1)-\frac{a c_2}{c_0+(i-2)c_2}Q_1(r_2+1)
-\frac{(i-a-1) c_2}{c_0+(i-2)c_2}Q_1(r_2-1)\bigg)
\nonumber\\
&&=\frac{1}{c_2}(n-2r_1+2)
((n-2r_2+2)(c_0-c_2)+4c_2a)\sum_{i=a+1}^{r_1}\left(\frac{1}{c_0+(i-2)c_2}-\frac{1}{c_0+(i-1)c_2)}\right)
\nonumber\\
&&=\frac{1}{c_2}(n-2r_1+2)
((n-2r_2+2)(c_0-c_2)+4c_2a)
\left(\frac{1}{c_0+(a-1)c_2}-\frac{1}{c_0+(r_1-1)c_2}\right)
\nonumber\\
&&=\frac{(n-2r_1+2)((c_0-c_2)(n-2r_2+2)+4c_2a)
(r_1-a)}
{(c_0+(a-1)c_2)(c_0+(r_1-1)c_2)}.
\label{equ:N4.56}
\end{eqnarray}
\begin{eqnarray}
&&\sum_{i=r_1+1}^{r_2+a}\frac{c_0+(i-2)c_2}{(c_0-c_2)(c_0+(i-1)c_2)}
\times
\nonumber\\
&&
\bigg(Q_1(r_1+1)-\frac{r_1c_2}{c_0+(i-2)c_2}Q_1(r_1-1)
-\frac{(i-r_1-1)c_2}{c_0+(i-2)c_2}Q_1(r_1+1)\bigg)\times
\nonumber\\
&&\bigg(Q_1(r_2-1)-\frac{a c_2}{c_0+(i-2)c_2}Q_1(r_2+1)
-\frac{(i-a-1)c_2}{c_0+(i-2)c_2}Q_1(r_2-1)\bigg)
\nonumber\\
&&=\frac{((c_0-c_2)(n-2r_2+2)+4c_2a)((c_0-c_2)(n-2)-2r_1(c_0+c_2))}{c_2(c_0-c_2)}
\times
\nonumber\\
&&\sum_{i=r_1+1}^{r_2+a}
\left(
\frac{1}{c_0+(i-2)c_2}-\frac{1}{c_0+(i-1)c_2}\right)
\nonumber\\
&&=\frac{(r_2-r_1+a)\big((c_0-c_2)(n-2r_2+2)
+4c_2 a\big)\big((n-2)(c_0-c_2)
-2r_1(c_0+c_2)\big)}{(c_0-c_2)(c_0+(r_1-1)c_2)
(c_0+(r_2-1)c_2+c_2a)}.
\nonumber\\
&&\label{equ:N4.57}
\end{eqnarray}
\begin{eqnarray}
&&\sum_{i=r_2+a+1}^{n}\frac{c_0+(i-2)c_2}{(c_0-c_2)(c_0+(i-1)c_2)}
\times
\nonumber\\
&&\bigg(Q_1(r_1+1)-\frac{r_1c_2}{c_0+(i-2)c_2}Q_1(r_1-1)
-\frac{(i-r_1-1)c_2}{c_0+(i-2)c_2}Q_1(r_1+1)\bigg)\times
\nonumber\\
&&\bigg(Q_1(r_2+1)-\frac{r_2c_2}{c_0+(i-2)c_2}Q_1(r_2-1)
-\frac{(i-r_2-1)c_2}{c_0+(i-2)c_2}Q_1(r_2+1)\bigg)\nonumber\\
&&=\frac{((n-2)(c_0-c_2)-2r_2(c_0+c_2))((n-2)(c_0-c_2)-2r_1(c_0+c_2))}{(c_0-c_2)c_2}\times
\nonumber\\
&&
\sum_{i=r_2+a+1}^n\left(\frac{1}{c_0+(i-2)c_2}-\frac{1}{c_0+(i-1)c_2}\right)
\nonumber\\
&&=\frac{\big((n-r_2)-a\big)\big((n-2)(c_0-c_2)-2r_2(c_0+c_2)\big)\big((n-2)(c_0-c_2)-2r_1(c_0+c_2)\big)}
{(c_0-c_2)\big(c_0+c_2(r_2-1)+c_2a\big)\big(c_0+(n-1)c_2\big)}.
\nonumber\\
&&\label{equ:N4.58}
\end{eqnarray}
Then (\ref{equ:N4.54}), (\ref{equ:N4.55}),
(\ref{equ:N4.56}), (\ref{equ:N4.57}) and (\ref{equ:N4.58}) implies 
{\small
\begin{eqnarray}&&
\frac{1}{(c_2-c_0)((n-1)c_2+c_0)}
\bigg\{16((n-1)c_2+c_0)a
-2r_1\bigg((n^2+4n-4)c_2\nonumber\\
&&
-(n^2-4n-4)c_0\bigg)-4r_1r_2\bigg((n-4)c_0-nc_2\bigg)
+(n-2)^2(n-2r_2)(c_2-c_0)\bigg\}
\nonumber
\\
&&+\frac{\bigg((n-2)(n-2r_1)d_0-c_0-(n-1)c_2\bigg)\bigg((n-2)(n-2r_2)d_0-c_0-(n-1)c_2\bigg)}
{\bigg((n-1)c_2+c_0\bigg)\bigg((W_{r_1}+W_{r_2})(c_0+(n-1)c_2)-nd_0^2\bigg)}\nonumber
\\
&&
=0.
\end{eqnarray}
}
\noindent
Then we have
$$\frac{(n-1)((n-r_2)r_1-na)}{
r_1(n-r_1)W_{r_1}+r_2(n-r_2)W_{r_2}}=0,
$$
and $a=\frac{(n-r_2)r_1}{n}$.
This implies (3).
\hfill\qed\\

\subsection{Regular semi-lattices and Geometric relative $t$-designs}
There is one more important property satisfied 
by geometric relative $t$-designs of
association schemes attached to regular semi-lattices. In \cite{Delsarte-1977}, Delsarte proved that 
if P-polynomial association scheme
has the property of regular semi-lattice and
also satisfies the Q-polynomial property, then
$(Y,w)$ is a
relative $t$-design with respect to a point $u_0$
if and only if $(Y,w)$ is a geometric relative $t$-design with respect to the regular semi-lattice.
Let $\mathfrak X=(X,\{R_i\}_{0\leq i\leq n})$ be
the P-polynomial scheme associated with a 
regular semi-lattice $\Lambda$.
Let $h$ be the hight function
of $\Lambda$ with $0\leq h(x)\leq n$. 
Let $\Lambda_j=\{x\in \Lambda\mid h(x)=j\}$
for $j=0,1,\ldots,n$.
Then $X$ is the top fiber 
$\Lambda_n=\{x\in \Lambda\mid h(x)=n\}$ of 
$\Lambda$.
Let $\chi\in \mathcal F(X)$.
Assume $\chi(x)\geq 0$ for any $x\in X$.
Let $j\in \{1,2,\ldots,n\}$ and we define
the following function $\lambda_{j,\chi}$
on $\Lambda_j$ by
\begin{eqnarray}
&&\lambda_{j,\chi}(z)=\sum_{x\in \Lambda_n\atop
{x\geq z}}\chi(x),~z\in \Lambda_j.
\end{eqnarray}
If the following condition satisfied, then 
$\chi$ is called
geometric relative $t$-design with respect to a point $u_0\in X$. 
For any integer satisfying $0\leq j\leq t$, there
exists a constant $\lambda_{u_0,j}$
and
\begin{eqnarray}
\lambda_{t,\chi}(z)=\lambda_{u_0,j}
\end{eqnarray}
holds for any $z\in \Lambda_t$ satisfying
$h(z\wedge u_0)=j$.
Now we consider the semi-lattice structure 
which gives $H(n,2)$.
Let $\Lambda=\{(x_1,\ldots,x_n)\mid x_i\in \{0,1\}
~\mbox{or}~ x_i=\cdot\}$. For $x=(x_1,\ldots,x_n),~y= (y_1,\ldots,y_n)\in \Lambda$, we deine $x\leq y$ if
$x_i=y_i$ or $x_i=\cdot$ for $1\leq i\leq n$.
We defined $h(x)=|\{i\mid x_i\in \{0,1\}\}|$. 
Then $\Lambda$ is a regular semi-lattice with
the hight function $h$. Clearly 
the top fiber is $\Lambda_n=F^n$ and
$\Lambda_n$ gives association scheme $H(n,2)$. 
Now we consider the geometric relative $t$-design with respect to $u_0=(0,0,\ldots,0)\in X$.
Let $z=(z_1,\ldots,z_n)\in \Lambda_t, ~h(z\wedge u_0)=t-j$.
Then $|\{i\mid z_i\in \{0,1\}\}|=t$,
$|\{i\mid z_i=0\}|=t-j$, $|\{i\mid z_i=1\}|=j$.
Then for $z\in \Lambda_t,~h(z\wedge u_0)=t-j$
we have
\begin{eqnarray}
\lambda_{u_0,t-j}=\lambda_{t,\chi}(z)=
\sum_{x\in \Lambda_n\atop {x\geq z}}\chi(x).
\end{eqnarray}
For any $u\in X_j$, we define
$\overline u=\{i\mid u_i=1\}=\{i_1,\ldots,i_j\}$.
Then there exists $z=(z_1,\ldots,z_n)\in\Lambda_t$ satisfying
$z_{i_1}=z_{i_2}=\cdots =z_{i_j}=1$ and
$|\{i\mid z_i=0\}|=t-j$.
Then $x\in \Lambda_n$ satisfies 
$x\geq z$ if and only if $\overline u\subset\overline x$.
Let $Y=\{y\in X\mid \chi(y)>0\}$ and $w(y)
=\chi(y)$ for $y\in Y$, then $(Y,w)$ is a 
relative $t$-design in the style of Definition \ref{def:1.1} (see more information in \cite{Delsarte-1976,Delsarte-1977,Bannai-B-2012}).
The argument given above implies the following proposition.
\begin{pro}
\label{pro:4.4}
Let $(Y,w)$ be a relative $t$-design in $H(n,2)$
with respect to $u_0=(0,0,\ldots,0)$. Then
for any $u\in X_j$
$$\sum_{y\in Y, \overline u\subset \overline y}w(y)
=\lambda_j
$$
is a constant depends only on $j$ ($0\leq j\leq t$).
Here $\overline x=\{i\mid x_i=1, 1\leq i\leq n\}$
for $x=(x_1,\ldots,x_n)\in X$.
\end{pro}
Please refer \cite{Delsarte-1976} for more information on regular semi-lattices.

\subsection{Proof of Theorem \ref{theo:2.2}}
Proposition \ref{pro:4.3} (2) and (3) imply
Theorem \ref{theo:2.2} (1).\\

\noindent
{\bf List of possible parameters for 
$\boldsymbol{ n\leq 30}$}\\
We first determined the parameter
set $\{n,r_1,r_2,N_{r_1},N_{r_2},\alpha_1,\alpha_2,\gamma,\frac{w_{r_2}}{w_{r_1}}\}$
according to the formula given in Proposition  
\ref{pro:4.3}.
If $(Y,w)$ is a relative $t$-design with respect to 
$u_0$, then $(Y,\mu w)$,  $(\mu w)(y)=\mu w(y)$ for $y\in Y$, is also a
relative $t$-design with respect to 
$u_0$ for any positive real number $\mu$.
Therefore in the following argument we assume $w_{r_1}=1$. 
We apply Proposition \ref{pro:4.4} and
determine $\lambda_1$ and $\lambda_2$.
For this purpose we count the elements in
the set
$\{(x,y)\mid x\in X_i, y\in Y\}$ for $i=1,2$.
Then we have
\begin{eqnarray}
w_{r_1}{r_1\choose i}N_{r_1}+w_{r_2}{r_2\choose i}N_{r_2}={n\choose i}\lambda_i, ~\mbox{for}~i=1,2.
\label{equ:4.63}
\end{eqnarray}
We note that if $w_{r_1}=w_{r_2}=1$,
then $\lambda_i=|\{y\in Y\mid \overline u\subset
\overline y\}|$ for any $u\in X_i$ for $i=1,2$ and
$\lambda_1,\lambda_2$ must be integers.
We list the feasible parameters 
$n,r_1,r_2,N_{r_1},N_{r_2},\alpha_1,\alpha_2,\gamma,w=w_{r_2}~(w_{r_1}=1),\lambda_1,\lambda_2$ for
a tight relative $2$-designs with respect to a point $u_0$ for $6\leq n\leq 30$
below. 
\begin{eqnarray}
\begin{array}{l || rr| rr| rrr| l |r r|c|}
n&r_1&r_2&N_{r_1}&N_{r_2}&\alpha_1&\alpha_2&\gamma&w&\lambda_1 &\lambda_2 &\\
\hline
6(1)&2&3&3&4&4&4&3&1&3&1&\circ\\
6(2)& 3&4&  4& 3& 4& 4& 3& 1&4&2&\circ\\
\hline
10(1) &  2&  5&  5& 6& 4& 6& 5& \frac{2}{3}&3&1&\times\\
10(2)&  4&  5&  5& 6& 6& 6& 5& 1&5&2&\circ\\
10(3)&  5&  6&  6& 5& 6& 6& 5& 1&6&3&\circ\\
10(4)&  5&  8&  6& 5& 6& 4& 5& \frac{3}{2}&9&6&\times\\
\hline
12(1)&  3&  4&  4& 9& 6& 6& 5& 1&4&1&\circ\\
12(2)    &  3&  8&  4& 9& 6& 6& 7& 1&7&4&\circ\\
12(3)&  4&  6&  9& 4& 6& 8& 6& \frac{3}{4}&\frac{9}{2}&
\frac{3}{2}&\times\\
12(4)&  4&  9&  9& 4& 6& 6& 7& 1&6&3&\circ\\
12(5)&  6&  8&  4& 9& 8& 6& 6&  \frac{4}{3}&10&6&\times\\
12(6)&  8&  9&  9& 4& 6& 6& 5& 1&9&6&\circ\\
\hline
14(1)&  2&  7&  7& 8& 4& 8& 7& \frac{1}{2}&3&1&\circ\\
14(2)&  6&  7&  7& 8& 8& 8& 7& 1&7&3&\circ\\
14(3)&  7&  8&  8& 7& 8& 8& 7& 1&8&4&\circ\\
14(4)&  7&  12&  8& 7& 8& 4& 7&  2&16&12&\circ\\
\hline
15(1)&  5&  6&  6& 10& 8& 8& 7& 1&6&2&\circ\\
15(2)&  5&  9&  6& 10& 8& 8& 8& 1&8&4&\circ\\
 15(3)&  6&10&10&   6& 8& 8& 8& 1&8&4&\circ\\
    15(4)&  9&10&10&   6& 8& 8& 7& 1&10&6&\circ\\
\hline
18(1)  &  2&  9&  9& 10&  4& 10&  9& \frac{2}{5}&3&1&\times\\
18(2)    &  8&  9&  9& 10& 10& 10& 9& 1&9&4&\circ\\
 18(3)&  9&10&10&   9& 10& 10& 9& 1&10&5&\circ\\
18(4)      &  9&16& 10&  9& 10&   4& 9& \frac{5}{2}&25&20&\times\\
\hline
20(1)& 4& 5  & 5 &16&   8&   8&   7&1&5&1&\circ
\\
20(2) & 4& 15& 5 &16&   8&   8& 13& 1&13&9&\circ\\
20(3)& 5& 8  &16&  5&   8& 12&   9&\frac{2}{3}&\frac{16}{3}&\frac{4}{3}&\times\\
 20(4)   & 5& 12&16&  5&   8& 12& 11&\frac{2}{3}&6&2&\times\\
 20(5)   & 5& 16&16&  5&   8&   8& 13&1&8&4&\circ\\
 20(6)& 8& 15&  5&16& 12&   8& 11&\frac{3}{2}&20&14&\times\\
20(7)   & 12& 15&  5&16& 12&   8& 9&\frac{3}{2}&21&15&\times\\
 20(8) & 15& 16&  16&5& 8&   8& 7&1&16&12&\circ\\
\hline
\end{array}
\nonumber
\end{eqnarray}
\begin{eqnarray}
\begin{array}{l || rr| rr| rrr| l| r r|c|}
n&r_1&r_2&N_{r_1}&N_{r_2}&\alpha_1&\alpha_2&\gamma&
w&\lambda_1&\lambda_2&\\
\hline
21(1)  &3& 7 &  7&15&   6& 10&   8& \frac{3}{5}&4&1&\times\\
21(2)    &3&14&  7&15&   6& 10& 13& \frac{3}{5}&7&4&\times\\
21(3) &6& 7 &  7&15& 10& 10&   9&1&7&2&\circ\\
21(4)& 6&14& 7&  15& 10& 10& 12& 1&12&7&\circ\\
 21(5)&7& 9 &15&  7& 10& 12& 10&\frac{5}{6}&\frac{15}{2}&\frac{5}{2}&\times\\
21(6)  &7&12&15&  7& 10& 12& 11&\frac{5}{6}&\frac{25}{3}&\frac{10}{3}&\times\\
21(7)    &7&15&15&  7& 10& 10& 12&1&10&5&\circ\\
21(8)    &7&18&15&  7& 10&   6& 13&\frac{5}{3}&15&10&\times\\
 21(9)   & 9&14&  7&15& 12& 10& 11&\frac{6}{5}&15&9&\times\\
21(10)  &12& 14&7& 15& 12& 10& 10&\frac{6}{5}&16&10 &\times\\
21(11)&14& 15& 15& 7& 10& 10& 9& 1
&15&10&\circ\\
21(12)&14& 18&15& 7& 10& 6& 8&  \frac{5}{3}&20&15&\times\\
\hline
22(1)&  2&11&11&12
&4&12   & 11&\frac{1}{3}&3&1&\circ\\
22(2)    &10&11&11&12& 12& 12& 11&1&11&5&\circ\\
  22(3)  &11&12&12&11& 12& 12& 11&1&12&6&\circ\\
   22(4) &11&20&12&11& 12&   4&11 &3&36&30&\circ\\
\hline
24(1)&8&9&9& 16& 12& 12& 11& 1&9&3&\circ\\
24(2)    &8&15&9&16& 12& 12& 13&1&13&7&\circ\\
  24(3)  &9&16&16&9& 12& 12& 13&1&12&6&\circ\\
24(4)& 15& 16& 16& 9& 12& 12& 11& 1&16&10
&\circ\\
\hline
26(1)&  2&13&13&14&   4& 14& 13& \frac{2}{7}&3&1&\times\\
26(2) &  6&13&13&14& 10& 14& 13& \frac{5}{7}&8&3&\times\\
 26(3)  &  8&13&13&14& 12& 14& 13& \frac{6}{7}&10&4&\times\\
  26(4)    &12&13&13&14& 14& 14& 13&1&13&6&\circ\\
 26(5)     &13&14&14&13& 14& 14& 13&1&14&7&\circ\\
    26(6)  &13&18&14&13& 14& 12& 13& \frac{7}{6}&\frac{35}{2}&\frac{21}{2}&\times\\
   26(7)  &13&20&14&13& 14& 10& 13& \frac{7}{5}&21&14&\times\\
26(8)         &13&24&14&13& 14&   4& 13& \frac{7}{2}&49&42&\times\\
\hline
 27(1)  &9&15&7&21& 14& 14& 14&1&14&7&\times\\
27(2)&12&18&21& 7& 14& 14& 14& 1
&14&7&\times\\
\hline
28(1)&  7&  8&  8& 21& 12& 12& 11& 1&8&2&\circ\\
28(2)&  7&20&  8& 21& 12& 12& 17& 1&17&11&\circ\\
28(3) &  8&14&21&   8& 12& 16& 14& \frac{3}{4}&9&3&\times\\
 28(4)&  8&21&21&   8& 12& 12& 17& 1&12&6&\circ\\    
28(5) &14&20&  8& 21& 16& 12& 14& \frac{4}{3}&24&16&\times\\
28(6)&20&21&21&8&12& 12& 11& 1&21&15&\circ\\
\hline
30(1)& 2& 15&15&16& 4&16&15&\frac{1}{4}&3&1&\circ\\
30(2) & 3& 10&10&21&6&14&11&\frac{3}{7}&4&1&\times\\
30(3)   & 3& 20&10&21&6&14&19&\frac{3}{7}&7&4&\times\\
 30(4)   & 5& 6&6&25&10&10&9&1&6&1&\circ\\
   30(5) & 5& 24&6&25&10&10&21&1&21&16&\circ\\
30(6)     & 6& 10&25&6&10&16&12&\frac{5}{8}&\frac{25}{4}
 &\frac{5}{4}&\times\\
30(7)   & 6& 15&25&6&10&18    &15&\frac{5}{9}&\frac{20}{3}&\frac{5}{3}&\times\\
30(8)     & 6& 20&25&6&10&16&18&\frac{5}{8} &\frac{15}{2}
  &\frac{5}{2}&\times\\
  \end{array}
\nonumber
\end{eqnarray}
\begin{eqnarray}
\begin{array}{l || rr| rr| rrr| l |r r|c|}
n&r_1&r_2&N_{r_1}&N_{r_2}&\alpha_1&\alpha_2&\gamma&
w&\lambda_1&\lambda_2&\\
\hline
  30(9)  & 6& 25&25&6&10&10&21&1&10&5&\circ\\
   30(10) & 9& 10&10&21&14&14&13&1&10&3&\circ\\
    30(11)& 9& 20&10&21&14&14&17&1&17&10&\circ\\
   30(12)  & 10& 12&21&10&14&16&14&\frac{7}{8}&\frac{21}{2}&\frac{7}{2}&\times\\
30(13)    & 10& 18&21&10&14&16&16&\frac{7}{8}&\frac{49}{4}&\frac{21}{4}&\times\\
   30(14) & 10& 21&21&10&14&14&17&1&14&7&\circ\\
30(15)&10&24&6&25& 16& 10& 18&  \frac{8}{5}&34&26&\times\\
30(16)  & 10& 27&21&10&14&6&19&\frac{7}{3}&28&21&\times\\
 30(17)  & 12& 20&10&21&16&14&16&\frac{8}{7}&20&12&\times\\
   30(18) & 14& 15&15&16&16&16&15&1&15&7&\circ\\
    30(19)& 15& 16&16&15&16&16&15&1&16&8&\circ\\
 30(20)    & 15& 24&6&25&18&10&15&\frac{9}{5}&39&30&\times\\
 30(21) & 15& 28&16&15&16&4&15&4&64&56&\circ\\
30(22)&18&20&10& 21&16& 14& 14&\frac{ 8}{7}&22&14&\times\\
30(23)&20&21&21&10&14&14&13&1&21&14&\circ\\
30(24)&20&24&6& 25&16&10&12&\frac{8}{5}&36&28&\times\\
30(25)& 20& 27& 21& 10& 14& 6& 11& \frac{ 7}{3}&35&28&\times\\
30(26)&24&25&25&6&10&10&9&1&25&20&\circ\\
\hline
\end{array}
\nonumber
\end{eqnarray}
The last column in the table given above, `` $\circ$ " indicates existence, `` $\times$ " indicates 
non-existence of the tight relative $2$-design with the corresponding parameters. 
For the cases with `` $\circ$ ", complete classification problem is still open.\\

\noindent
{\bf Constructions}\\
First we give two kind of construction theorem.
First one is the construction by Hadamard 
matrices.\\
\noindent
Let $m\equiv -1~(\mbox{mod}~4)$,
and $n=2m$. Suppose there is an Hadamard matrix $H_{m+1}$ of size
$(m+1)\times (m+1)$. Let $h_1,h_2,\ldots,h_{m+1}$ be the row vectors of $H_{m+1}$. We may assume that each vector $h_j$ is of the following form by normalization, i.e.,
$h_j=(+,a_{j,1},a_{j,2},\ldots, a_{j,m})$ with
$a_{j,\nu}\in \{+,-\},~ 1\leq \nu\leq m$.   
First we define $Y_2\subset X_2$ in the following way.
\begin{eqnarray}&&Y_2=\{(y_{i,1},y_{i,2},\ldots,
y_{i,2\nu-1},y_{i,2\nu},\ldots,
y_{i,2m-1},y_{i,2m})\in X_2\mid
\nonumber\\
&&
 \quad1\leq i\leq m,~
y_{i,2i-1}=y_{i,2i}=1, y_{i,\nu}=0, \nu\neq 2i-1,~2i
\}.\nonumber
\end{eqnarray}
Then $|Y_2|=m=\frac{n}{2}$. Next we define
$Y_{m}$ in $X_m$ using $m+1$ row vectors
$h_1,\ldots,h_{m+1}$ of $H_{m+1}$.
For each $h_j$ ($1\leq j\leq m+1$), we define 
$$y(h_j)=(y_{j,1},
y_{j,2},\ldots,y_{j,2\nu+1},y_{j,2\nu+2},\ldots,
y_{j,2m-1},y_{j,2m})\in X_m$$
as follws: the $(2\nu-1)$-th and $2\nu$-th 
entries 
$y_{j,2\nu-1},y_{j,2\nu}$ of $y(h_j)$ are given by
$$(y_{j,2\nu-1},y_{j,2\nu})=
\left\{\begin{array}{ll}(1,0)&\mbox{if $a_{j,\nu}=+$}\\
(0,1)&\mbox{if $a_{j,\nu}=-$}
\end{array}
\right.
$$
for $\nu=1,\ldots,m$. Let $Y_m=\{y(h_1),\ldots,
y(h_{m+1})\}$.
Thens we have $|Y_m|=m+1$ and $Y=Y_2\cup Y_m$
satisfies the conditions of relative $2$-design
with respect to $u_0=(0,0,\ldots,0)$ in $H(n,2)$.
We can also easily check that the set
$Y'=\{(1,1,\ldots,1)-y\mid y\in Y\}$ which is the complement of $Y$ also is a relative 2-design
with respect to $u_0$ with 
$w'(y')=\frac{1}{w((1,1,\ldots,1)-y')}$ for $y'\in Y'$. This completes the proof of Theorem 
\ref{theo:2.2} (3).\\

Next one is the constructions from symmetric designs.
The following proposition is known.
\begin{pro}[Woodall \cite{Woodall-1970}]
\label{pro:4.5}
 Let $(V,\mathcal B)$ be a symmetric design $2$-$(n+1,k,\lambda)$ design. Let the point set
 $V=\{0,1,2,\ldots,n\}$.
\begin{enumerate}
\item Let $Y_{k-1}=\{y\in F^n\mid \overline y=
B\backslash \{0\}, B\in \mathcal B,  0\in B\}$
and $Y_{k}=\{y\in F^n\mid \overline y=B, B\in \mathcal B, 0\not\in B\}$. Then $Y=Y_{k-1}\cup Y_k$
is a tight relative $2$-design of $H(n,2)$ with respect to $u_0=(0,0,\ldots,0)$.
\item Let $2k\not=n+1$. Let
$Y_{n-k+1}=\{y\in F^n\mid \overline y= V\backslash B, B\in \mathcal B, 0\in B\}$ 
and $Y_k=\{y\in F^n\mid \overline y= B, B\in
\mathcal B, 0\not\in B\}$.
Then $Y=Y_{n-k+1}\cup Y_k$ is a tight relative $2$-design withrespect to $u_0$.
\end{enumerate}
\end{pro}
It is known that the complement of a symmetric design is also a symmetric design.
Therefore
using Proposition \ref{pro:4.5}, we can
construct tight relative $2$-design of $H(n,2)$ with respect to $u_0=(0,0,\ldots,0)$
for each set of parameters in the table
satisfying $w=1$ except for $n=27$.

\begin{rem}
\begin{enumerate}
\item
The tight relative $2$-designs with the parameters 
14(1), 14(4), 22(1), 22(4), 30(1) and 30(21)
are constructed by using Hadamard matrices
according to the method given above (Theorem \ref{theo:2.2} (3)).  
\item For $n=27$ the parameters do not correspond to
symmetric $2$-$(n+1,k,\lambda)$ designs.
\end{enumerate}
\end{rem}

\noindent
{\bf Non-existence}\\
In the following we prove that for each 
set of parameters 
with `` $\times$ " in the last column, 
tight relative $2$-design does not exist.
\begin{pro}
\label{pro:4.7N}
\begin{enumerate}
\item  Let $u\in X_1$ and 
$\lambda^{(i)}_1(u)=|\{y\in Y_{r_i}\mid
\overline  u\subset \overline y \}|$ for $i=1,2$.
Then $\lambda^{(i)}_1(u)$ does not depend on the choice of $u\in X_1$ and given by the following formulas.
\begin{enumerate}
\item 
\begin{eqnarray}&&\lambda^{(1)}_1=\lambda^{(1)}_1(u)
=\frac{ (r_2-1)\lambda_1- (n-1)\lambda_2}
{(r_2-r_1)w_1},
\label{equ:4.64}\\
&&\lambda^{(2)}_1=\lambda^{(2)}_1(u)
= \frac{(n-1)\lambda_2 -(r_1-1)\lambda_1}
{(r_2-r_1)w_2}.\label{equ:4.65}
\end{eqnarray}
\item The following holds.
\begin{eqnarray}&&
\sum_{i=1}^2\frac{N_{r_i}w_{r_i}}{|X_{r_i}|}
\left({n-1\choose r_i-1}Q_1(r_i-1)
+{n-1\choose r_i}Q_1(r_i+1)\right)
\nonumber\\
&&=\sum_{i=1}^2w_{r_i}
\bigg(\lambda^{(i)}_1Q_1(r_i-1)
+(N_{r_i}-\lambda^{(i)}_1)Q_1(r_i+1)\bigg).
\label{equ:4.66}
\end{eqnarray}
\end{enumerate}
\item Let $u\in X_2$, and let 
$\lambda^{(i)}_2(u)=|\{y\in Y_{r_i}\mid \overline u\subset \overline y\}$ and
$\lambda^{(i)}_{2,C}(u)=|\{y\in Y_{r_i}\mid \overline u\cap \overline y=\emptyset \}$ for $i=1,2$.
Then the following holds.
\begin{eqnarray}&&\sum_{i=1}^2
\frac{N_{r_i}w_{r_i}}{|X_{r_i}|}
\left({n-2\choose r_i-2}Q_2(r_i-2)+{n-2\choose r_i}Q_2(r_i+2)+2{n-2,\choose r_i-1}Q_2(r_i)\right)
\nonumber\\
&&=\sum_{i=1}^2w_{r_i}
\bigg(\lambda^{(i)}_2(u)Q_2(r_i-2)+\lambda^{(i)}_{2,C}(u)Q_2(r_i+2)\nonumber\\
&&\hskip4cm
+\big(N_{r_i}-\lambda^{(i)}_2(u)-\lambda^{(i)}_{2,C}(u)\big)Q_2(r_i)
\bigg).
\label{equ:4.67}
\end{eqnarray}
\end{enumerate}
\end{pro}
{\bf Proof}
(1) (a) Let $u\in X_1$ be fixed arbitrarily and consider the following
sum.
\begin{eqnarray}
\sum_{y\in Y,~\overline u\subset \overline y
\atop
{\{\overline x,\overline u\}\subset \overline y},~
x\in X_1,~x\not= u }w(y)
=\sum_{x\in X_1,\atop{x\not= u}}
\sum_{y\in Y,\atop{
\{\overline x,\overline u\}\subset \overline y }}w(y)
=\sum_{x\in X_1,\atop{x\not= u}}
\lambda_2=(n-1)\lambda_2.
\label{equ:4.68}
\end{eqnarray}
The left side of (\ref{equ:4.68}) has the following
reformation.
\begin{eqnarray}&&
\sum_{y\in Y,~\overline u\subset \overline y
\atop
{\{\overline x,\overline u\}\subset \overline y},~
x\in X_1,~x\not= u }w(y)
=\sum_{i=1}^2\sum_{y\in Y_{r_i},
\atop{\overline u\subset \overline y}}
\sum_{x\in X_1, x\not=u,
\atop{\overline x\subset \overline y}}w(y)
=\sum_{i=1}^2(r_i-1)w_{r_i}
|\{y\in Y_{r_i}\mid \overline u\subset \overline y\}|
\nonumber\\
&&=\sum_{i=1}^2(r_i-1)w_{r_i}\lambda^{(i)}_1(u).
\end{eqnarray}
Therefore we must have
\begin{eqnarray}
\sum_{i=1}^2(r_i-1)w_{r_i}\lambda^{(i)}_1(u)=(n-1)\lambda_2.
\label{equ:4.70}
\end{eqnarray}
On the other hand Proposition \ref{pro:4.4}
implies 
\begin{eqnarray}\lambda^{(1)}_1(u)w_{r_1}+\lambda^{(2)}_1(u)w_{r_2}=\lambda_1.
\label{equ:4.71}
\end{eqnarray}
Since the coefficient matrix of 
equations (\ref{equ:4.70})
and (\ref{equ:4.71}) with variable $\lambda^{(1)}_1(u),\lambda^{(2)}_1(u)$ is nonsingular 
$\lambda^{(1)}_1(u),\lambda^{(2)}_1(u)$
are determined by the
formulas in (1) (a).
\\
(1) (b) and 
(2): The 
equation (\ref{equ:1.4}) for  $\phi_u^{(1)},~u\in X_1$ and 
$\phi_u^{(2)},~u\in X_2$ imply equation
(\ref{equ:4.66}) and (\ref{equ:4.67}) respectively.
\hfill\qed \\

\noindent
Proposition 2.2 (2) in \cite{Li-B-B-2013} implies that
if (\ref{equ:4.66}) and  (\ref{equ:4.67}) 
are satisfied for each $u\in X_1$ and 
$u\in X_2$ respectively, then
$(Y_{r_1}\cup Y_{r_2} ,w)$ is a relative 2-design. 
On the other hand Proposition \ref{pro:4.4}
implies 
\begin{eqnarray}\lambda^{(1)}_2(u)w_{r_1}+\lambda^{(2)}_2(u)w_{r_2}=\lambda_2
\label{equ:4.73N}
\end{eqnarray}
for any $u\in X_2$.
In the following we will show that
for each case marked `` $\times$ ", there is no set
of integers $\{\lambda^{(1)}_1, \lambda^{(2)}_1\}\cup
\{\lambda^{(1)}_2(u), \lambda^{(1)}_{2,C}(u), \lambda^{(2)}_2(u),
 \lambda^{(2)}_{2,C}(u)
 \mid u\in X_2\}$ satisfying
(\ref{equ:4.66}), (\ref{equ:4.67}) and (\ref{equ:4.73N}). \\

\noindent
{\bf $\boldsymbol{n=10}$: non-existence for 10(1) and 10(4)}\\
$\bullet$ 10(1): Equation (\ref{equ:4.67}) implies
$$2\lambda^{(1)}_2(u)+2-\lambda^{(1)}_{2,C}(u)
+\frac{1}{3}\lambda^{(2)}_2(u)+\frac{1}{3}\lambda^{(2)}_{2,C}(u) = 0.$$
On the other hand 
(\ref{equ:4.73N}) implies 
 $ \lambda^{(1)}_2(u)+\frac{2}{3}\lambda^{(2)}_2(u)=\lambda_2=1$.
Since $0\leq \lambda^{(1)}_2,\lambda^{(1)}_{2,C}\leq 5$,
$\lambda^{(1)}_2+\lambda^{(1)}_{2,C}\leq 5$,
 $0\leq \lambda^{(2)}_2,\lambda^{(2)}_{2,C}\leq 6$, and
 $\lambda^{(2)}_2+\lambda^{(2)}_{2,C}\leq 6$,
only solution for these equations is
$\lambda^{(1)}_2(u)=1, \lambda^{(1)}_{2,C}(u)=4, \lambda^{(2)}_2(u)=\lambda^{(2)}_{2,C}(u)=0$. This contradict $r_2=5$.\\
$\bullet$ $10(4)$: We have 2 solutions 
$\lambda^{(1)}_2(u)=\lambda^{(1)}_{2,C}(u)=0, \lambda^{(2)}_2(u)=4,\lambda^{(2)}_{2,C}(u)=1$ and
$\lambda^{(1)}_2(u)=0, \lambda^{(1)}_{2,C}(u)=6, \lambda^{(2)}_2(u)=4, \lambda^{(2)}_{2,C}(u)=0$. This contradict $r_1=2$.\\

\noindent{\bf $\boldsymbol{n=12}$: non existence for 12(3) and 12(5)}\\
$\bullet$ $12(3)$: We have
$24\lambda^{(1)}_2(u)-\frac{40}{11}-8\lambda^{(1)}_{2,C}(u)
+8\lambda^{(2)}_2(u)+8\lambda^{(1)}_{2,C}(u) = 0$ and
$25\lambda^{(1)}_2(u)+\frac{3}{4}\lambda^{(2)}_2(u) =\lambda_2 =\frac{3}{2}$. Then
$\lambda^{(1)}_2(u) = 0, \lambda^{(1)}_{2,C}(u)= 3, \lambda^{(2)}_2(u)= 2, \lambda^{(1)}_{2,C}(u)= 2$ is the unique solution of these equations. This contradicts $r_1=4$.\\
$\bullet 12(5)$: 
$r_1=6,~r_2=8$, $N_{r_1}=4,~N_{r_2}=9$.
We obtain 
$\lambda^{(1)}_1=2,
\lambda^{(2)}_1=6$ and 
$\lambda^{(1)}_2(u)=2,
\lambda^{(2)}_2(u)=3$
 for any $u\in X_2$.
 Let $Y_{r_1}=\{y_1,\ldots,y_4\}$.
Since $r_1=6$
and $\alpha_1=8$, 
we have $|\overline y_i\cap \overline y_j|=2$
for $i\neq j$. Hence
we may assume
$\overline y_1=\{1,2,3,4,5,6\}$
and $\overline y_2=\{1,2,7,8,9,10\}$.
Since $\lambda^{(1)}_1=2$, we must have
$1,2\not\in  \overline y_3, \overline y_4$.
Then $\{1,3\}\not\subset \overline y_2,~\overline y_3,~\overline y_4$.
This contradicts $\lambda^{(1)}_2(\{1,3\})=2$.
\hfill\qed\\

\noindent{\bf $\boldsymbol{n=18}$: non existence for 18(1) and 18(4).}\\
$\bullet $ $18(1) $:
Similar computation shows that there is unique solution
$\lambda^{(1)}_2(u) = 1, \lambda^{(1)}_{2,C}(u)  = 8, \lambda^{(2)}_2(u) =  \lambda^{(2)}_{2,C}(u)  = 0$. This contradicts $r_2=9$.\\
$\bullet $ 18(4):
Similar computation shows that there is unique solution
$\lambda^{(1)}_2(u) =  \lambda^{(1)}_{2,C}(u)  = 0, \lambda^{(2)}_2(u) =8,
  \lambda^{(2)}_{2,C}(u)  =1$. This contradicts $r_1=9$.\\

\noindent{\bf $\boldsymbol{n=20}$: non existence for 20(3), 20(4), 20(6),  20(7).}\\
$\bullet $ $20(3)$:
Similar computation shows that there is unique solution
$\lambda^{(1)}_2(u) = 0, \lambda^{(1)}_{2,C}(u)  = 8, \lambda^{(2)}_2(u) =2,
  \lambda^{(2)}_{2,C}(u)  =3$. This contradicts $r_1=5$.\\
$\bullet$ $20(4)$:
Similar computation shows that there 
$\lambda^{(1)}_2(u) = 0$ or $\lambda^{(1)}_2(u) = 2$.
On the other hand $\alpha_1=8$ implies $|y_1\cap y_2|=1$ for any $y_1,y_2\in Y_{r_1}$. 
Then we must have $\lambda^{(1)}_2(u) = 0$ and contradicts $r_1=5$. \\
$\bullet $ $20(6)$:  
$r_1=8$, $r_2=15$, $N_1=5$, $N_2=16$.
We have $\lambda^{(1)}_1=2$ and 
$\lambda^{(2)}_1=12$, and
$\lambda^{(1)}_2(u)= 2,~ 
\lambda^{(1)}_{2,C}(u)= 3,~ 
\lambda^{(2)}_2(u)= 8,~ 
\lambda^{(2)}_{2,C}(u)= 0$.
Let $Y_8=\{y_1,\ldots,y_5\}$.
Since $\alpha_1=12$, we must have 
$|\overline y_i\cap \overline y_j|=2$
for any distinct $y_i, y_j\in Y_{r_1}$.
Then we may assume $\overline y_1=\{1,2,3,4,5,6,7,8\}$ and $\overline y_2=\{1,2,9,10,11,12,13,14\}$.\\
Since $\lambda^{(1)}_1=2$, we must have
$1\not\in \overline y_j$ for $j=3,4,5$.
Hence $\lambda^{(1)}_2(\{1,3\})= 2$
implies $\{1,3\}\subset \overline y_1, \overline y_2$. But this is impossible.\\
$\bullet $ $20(7)$:  $r_1=12$, $N_{r_1}=5$.
In this case we have
the following solutions.
$\lambda^{(1)}_1 = 3$, $\lambda^{(2)}_1 = 12$, and 
$\lambda^{(1)}_2(u) = 0$, or $3$ for any $u\in X_2$. 
Let $Y_{r_1}=\{y_1,y_2,\ldots,y_5\}$.
Since 
$\alpha_1=12$, we must have $|\overline y_i\cap\overline y_j|=6$ for any distinct $y_i$, $y_j\in Y_{r_1}$.
We may assume $\overline y_1=\{1,2,3,4,5,6,7,8,9,10,11,12\}$
and $\overline y_2=\{1,2,3,4,5,6,
13,14,15,16,17,18\}$.\\
Since $\lambda^{(1)}_1 = 3$, we may assume
$1\in \overline y_3$, then we must have
$1\not\in \overline y_4, \overline y_5$.
On the other hand, $\lambda^{(1)}_2(u) = 0$, or $3$ and $\{1,7\}\subset \overline y_1$
implies $\lambda^{(1)}_2(\{1,7\})=3$.
Since $1\not\in \overline y_4, \overline y_5$,
we must have $\{1,7\}\subset \overline y_1,\overline y_2,\overline y_3$.
This is impossible.\\

\noindent{\bf $\boldsymbol{n=21}$: non existence for 21(1), 21(2), 21(5), 21(6), 21(8), 21(9),21(10),21(12).}\\
$\bullet $ $21(1)$:
 We have $\lambda^{(2)}_2(u)=0$ for any $u\in X_2$.
 This contradicts $r_2=7$.\\
 \noindent
$\bullet $ $21(2)$: $r_1=3$ and $N_{r_1}=7$.
 Similar computation implies there is unique solution
 $\lambda^{(1)}_2(u)  = 1, \lambda^{(1)}_{2,C}(u) = 6, \lambda^{(2)}_2(u)  = 5, \lambda^{(2)}_{2,C}(u)  = 0$.
 However in this case we have $\alpha_1=6$ implies 
 $|\overline{y_1}\cap \overline{y_2}|=0$ 
 for any distinct $y_1,y_2\in Y_{r_1}$.
 This shows that there exists $u\in X_2$
 satisfying
 $\lambda^{(1)}_2(u)  =0$. This is a contradiction.
 \\
 $\bullet $ $21(5)$:
There is a unique solution 
$\lambda^{(1)}_2(u)  = 0, \lambda^{(1)}_{2,C}(u)  = 5, \lambda^{(2)}_2(u)  = 3, \lambda^{(2)}_{2,C}(u)  = 4$. This contradicts $r_1=7$.\\
$\bullet $ $21(6)$.
There is a unique solution
$\lambda^{(1)}_2(u) = 0, \lambda^{(1)}_{2,C}(u)  = 5,\lambda^{(2)}_2(u) = 4, \lambda^{(2)}_{2,C}(u)= 3$. This contradicts $r_1=7$.\\
$\bullet $ $21(8)$:
 There is a unique solution
 $\lambda^{(1)}_2(u)= 0,\lambda^{(1)}_{2,C}(u)= 5, \lambda^{(2)}_2(u)= 6, \lambda^{(2)}_{2,C}(u) = 1$. This contradicts $r_1=7$.\\
 $\bullet $ $21(9)$: $r_1=9$ and $N_{r_1}=7$. We have 
$\lambda^{(1)}_1 = 3, \lambda^{(2)}_1 = 10$ and 
$\lambda^{(1)}_2(u) = 3,  \lambda^{(2)}_2(u)  = 5$
for any $u\in X_2$. 
Let $Y_{r_1}=\{y_1,y_2,\ldots, y_7\}$.
Since $\alpha_1=12$,
$|\overline y_i\cap \overline y_j|=3$ for any 
distinct $y_i,y_j\in Y_{r_1}$.
We may assume $\overline y_1
=\{1,2,3,4,5,6,7,8,9\}$
and $\overline y_2=\{1,2,3,10,11,12,13,14,15\}$.
Since $\lambda^{(1)}_1 = 3$, we may assume 
$1\in\overline y_3$ and $1\not \in \overline y_j,~j=4,5,6,7$.
Hence $\lambda^{(1)}_2(\{1,4\}) = 3$ implies
$\{1,4\}\subset  \overline y_j$, for $j=1,2,3$.
This is impossible.\\
$\bullet $ $21(10)$: $r_1=12$, $N_{r_1}=7$. We have 
$\lambda^{(1)}_1= 4, \lambda^{(2)}_1 = 10$
and 
$\lambda^{(1)}_2(u) = 4,  \lambda^{(2)}_2(u)  = 5$
for any $u\in X_2$. 
Let $Y_{r_1}=\{y_i\mid 1\leq i\leq 7\}$.
Since $r_1=12$ and
$\alpha_1=12$, 
$|\overline y_i\cap \overline y_j|=6$ for any 
distinct $y_i,y_j\in Y_{r_1}$.
We may assume 
$\overline y_1
=\{i\mid 1\leq i\leq 12\}$
and $\overline y_2=\{1,2,3,4,5,6,13,14,15,16,17,18\}$.
Since $\lambda^{(1)}_1=4$ 
we may assume $1\in  \overline y_3, 
\overline y_4$ and $1\not\in  \overline y_j$, $j= 5,6,7$. 
Then $\lambda^{(1)}_2(\{1,7\})=4$
and $1\not\in  \overline y_j$, $j= 5,6,7$
imply $\{1,7\}\subset \overline y_j$, $j=1,2,3,4$.
This is a contradiction.\\
$\bullet $ $21(12)$: $r_2=18$, $N_{r_2}=7$. We have 
$\lambda^{(1)}_1 = 10, \lambda^{(2)}_1 = 6$ 
and $\lambda^{(1)}_2(u) = 5, \lambda^{(2)}_2(u)  = 6$
for any $u\in X_2$. 
Let $Y_{r_2}=\{y_i\mid 1\leq i\leq 7\}$.
Since $r_2=18$ and $\alpha_2=6$,
$|\overline y_i\cap \overline y_j|=15$ for any 
distinct $y_i,y_j\in Y_{r_2}$.
We may assume 
$\overline y_1=\{i\mid 1\leq i\leq 18\}$
and $\overline y_2=\{1,2,\ldots,15,19,20,21\}$.
Since $\lambda^{(2)}_1=6$, we may assume
$1\in \overline y_j$, for $j=1,2,3,4,5,6$ and 
$1\not\in \overline y_7$.
Then it is imposible to have $\lambda^{(2)}_2(\{1,16\})=6$.
\\

\noindent{\bf $\boldsymbol{n=26}$: nonexistence for 26(i), $\boldsymbol{i\neq 4, 5}$}\\
$\bullet $ $26(1)$, $26(2)$ and $26(3)$: We have $\lambda^{(2)}_2(u)= 0$.
This contradicts $r_2=13$.\\
$\bullet $ $26(6)$, $26(7)$ and $26(8)$: We have $\lambda^{(1)}_2(u)= 0$. This contradicts $r_1=13$.\\

\noindent{\bf $\boldsymbol {n=27}$: nonexistence for 27(1),  27(2)}\\
$\bullet$ $27(1)$:  
$\lambda^{(1)}_1=\frac{7}{3}$ and
$\lambda^{(2)}_1=\frac{35}{3}$.
This is a contradiction.
\\
$\bullet$ $27(2)$:  
$\lambda^{(1)}_1=\frac{28}{3}$ and
$\lambda^{(2)}_1=\frac{14}{3}$.
This is a contradiction.
\hfill\qed\\ 

\noindent{\bf $\boldsymbol {n=28}$: nonexistence for 28(3) and 28(5)}\\
$\bullet $ $28(3)$: $r_2=14$, $N_{r_2}=8$. 
We have 
$\lambda^{(1)}_1= 6, \lambda^{(2)}_1=4$
and
$\lambda^{(1)}_2(u) = 0,  \lambda^{(2)}_2(u)  = 4$;
or $\lambda^{(1)}_2(u) = 3,  \lambda^{(2)}_2(u)  = 0$
for any $u\in X_2$. 
Let $Y_{r_2}=\{y_i\mid 1\leq i\leq 8\}$.
Then $r_2=14$ and $\alpha_2=16$
implies 
$|\overline y_i\cap \overline y_j|=6$ for any 
distinct $y_i,y_j\in Y_{r_2}$.
We may assume 
$\overline y_1=\{i\mid 1\leq i\leq 14\}$
and $\overline y_2=\{1,2,3,4,5,6,15,16,17,18,19,20,21,22\}$.
Since $\lambda^{(2)}_2(\{1,2\}) \geq 1$ 
we must have $\lambda^{(2)}_2(\{1,2\})=4$.
Therefore we may assume $\{1,2\}\subset \overline y_3, \overline y_4$.
Then $\lambda^{(2)}_1=4$ implies
$1\not\in \overline y_j$ for $j=5,6,7,8$.
$\{1,7\}\subset \overline y_7$ implies
$\lambda^{(2)}_2(\{1,7\})=4$.
Therefore we must have $\{1,7\}\subset
\overline y_j$ for $j=1,2,3,4$. But 
this is impossible.\\
 $\bullet $ $28(5)$: $r_1=14$ and $N_{r_1}=8$.
 We have 
$\lambda^{(1)}_1= 4, \lambda^{(2)}_1=15$ and 
$\lambda^{(1)}_2(u) = 0,  \lambda^{(2)}_2(u)  = 12$;
or $\lambda^{(1)}_2(u) = 4,  
\lambda^{(2)}_2(u)  =9$ for any $u\in X_2$. 
Let $Y_{r_1}=\{y_i\mid 1\leq i\leq 8\}$.
Then $r_1=14$ and $\alpha_1=16$
implies 
$|\overline y_i\cap \overline y_j|=6$ for any 
distinct $y_i,y_j\in Y_{r_1}$.
We may assume $\overline y_1=\{1,2,\ldots,14\}$ and $\overline y_2=\{1,2,\ldots,6,15,16,17,18,19,20,21,22\}$.
Since $\lambda^{(1)}_1= 4$, we may assume
$1\in \overline y_3, \overline y_4$ and
$1\not\in \overline y_j$, $j=5,4,3,8$.
Since $\{1,14\}\subset\overline y_1$, we must
have $\lambda^{(1)}_2(\{1,14\}) =4$.
However since $14\not\in \overline y_2$
this is impossible.\\

\noindent{\bf $\boldsymbol {n=30}$:  non existences for `` $\times$ "}\\
$\bullet $ $30(2)$: We have $\lambda^{(2)}_2(u)=0$
for any $u\in X_2$.
This contradicts $r_2=10$.\\
$\bullet $ $30(3)$: We have $\lambda^{(1)}_2(u)=1$.
Since $r_1=3$ and $\alpha_1=6$,
there exists $u\in X_2$ with $\lambda^{(1)}_2(u)=0$.
This is a contradiction.\\
$\bullet $ $30(6)$, $30(7)$ and $30(8)$:  We have $\lambda^{(1)}_2(u)=0$ for any $u\in X_2$.
This contradicts $r_1=6$.\\ 
$\bullet $ $30(12)$, $30(13)$ and $30(16)$:  We have $\lambda^{(1)}_2(u)=0$ for any $u\in X_2$.
This contradicts $r_1=10$.\\ 
$\bullet $ $30(15)$: $r_1=10$ 
and $N_{r_1}=6$.
We have $\lambda^{(1)}_1 = 2, \lambda^{(2)}_1= 20$ and $\lambda^{(1)}_2(u)= 2,
\lambda^{(2)}_2(u)= 15$ for any $u\in X_2$.
Let $Y_{r_1}=\{y_1,\ldots,y_6\}$.
Then $r_1=6$ and $\alpha_1=16$,
implies $|\overline y_i\cap \overline y_j|=8$ for any 
distinct $y_i,y_j\in Y_{r_1}$. 
We may assume 
$\overline y_1=\{1,2,3,4,5,6,7,8,9,10\}$,
$\overline y_2=\{1,2,3,4,5,6,7,8,11,12\}$.
Since $\lambda^{(1)}_1=2$, we must have
$1,2\not\in\overline y_i$, for $i=3,4,5,6$.
Then it is impossible to have 
 $\lambda^{(1)}_2(\{1,9\})= 2$.\\
 $\bullet $ $30(17)$: $r_1=12$ and $N_{r_1}=10$.
We have $\lambda^{(1)}_1 = 4, \lambda^{(2)}_1=14$ and $\lambda^{(1)}_2(u)=4,
\lambda^{(2)}_2(u)= 7$ for any $u\in X_2$.
Let $Y_{r_1}=\{y_1,\ldots,y_{10}\}$.
Then $r_1=12$ and $\alpha_1=16$,
implies $|\overline y_i\cap \overline y_j|=4$ for any 
distinct $y_i,y_j\in Y_{r_1}$. 
We may assume 
$\overline y_1=\{1,2,3,4,5,6,7,8,9,10,11,12\}$,
$\overline y_2=\{1,2,3,4,13,14,15,16,17,18,19,20\}$.
Since $\lambda^{(1)}_2=4$ and $\lambda^{(1)}_1=4$, we may assume
$\{1,2\}\subset \overline y_j$ for $j=3,4$
and $1,2\not\in \overline y_j$ for $j=5,6,7,8,9,10$.
On the other hand $\lambda^{(1)}_2(\{1,5\})=4$.
Hence we must have $\{1,5\}\subset
\overline y_j$ for $j=1,2,3,4$. But this is impossible.
\\
 $\bullet $ $30(20)$: $r_1=15$ and $N_{r_1}=6$.
 We have $\lambda^{(1)}_1= 3, \lambda^{(2)}_1=20$ and $\lambda^{(1)}_2(u)=3,
\lambda^{(2)}_2(u)= 15$ for any $u\in X_2$.
Let $Y_{r_1}=\{y_1,\ldots,y_6\}$.
Then $r_1=15$ and $\alpha_1=18$,
implies $|\overline y_i\cap \overline y_j|=4$ for any 
distinct $y_i,y_j\in Y_{r_1}$. 
We may assume 
$\overline y_1=\{i\mid 1\leq i\leq 15\}$,
$\overline y_2=\{1,2,3,4,16,17,18,19,20,21,22,23,24,\\
25, 26\}$.
Since $\lambda^{(1)}_2=3$ and $\lambda^{(1)}_1=3$, we may assume
$\{1,2\}\subset \overline y_3$ 
and $1,2\not\in \overline y_j$ for $j=4,5,6$.
Then it is impossible to have $\lambda^{(1)}_2(\{1,5\})=3$.\\
 $\bullet $ $30(22)$: $r_1=18$ and $N_{r_1}=10$.
 We have $\lambda^{(1)}_1 = 6, 
 \lambda^{(2)}_1=14$ and $\lambda^{(1)}_2(u)=6,
\lambda^{(2)}_2(u)= 7$ for any $u\in X_2$.
Let $Y_{r_1}=\{y_1,\ldots,y_{10}\}$.
Then $r_1=18$ and $\alpha_1=16$,
implies $|\overline y_i\cap \overline y_j|=10$ for any 
distinct $y_i,y_j\in Y_{r_1}$. 
We may assume 
$\overline y_1=\{i\mid 1\leq i\leq 18\}$,
$\overline y_2=\{1,2,\ldots,10,19,\ldots,26\}$.
Since $\lambda^{(1)}_2=\lambda^{(1)}_1=6$, we may assume
$\{1,2\}\subset \overline y_i$ for $i=3,4,5,6$ and
$1,2\not\in \overline y_i$ for $i=7,8,9,10$. 
Then we must have $j\in \overline y_i$
for $j=1,2,\ldots,10$ and $i=1,2,\ldots,6$.
Then it is impossible to have 
$\lambda^{(1)}_2(\{1,11\})=6$.\\
 $\bullet $ $30(24)$: $r_1=20$ and $N_{r_1}=6$.
 We have $\lambda^{(1)}_1 = 4,
  \lambda^{(2)}_1=20$ and $\lambda^{(1)}_2(u)=4, \lambda^{(2)}_2(u)= 15$ for any $u\in X_2$.
Let $Y_{r_1}=\{y_1,\ldots,y_{6}\}$.
Then $r_1=20$ and $\alpha_1=16$,
implies $|\overline y_i\cap \overline y_j|=12$ for any 
distinct $y_i,y_j\in Y_{r_1}$. 
We may assume 
$\overline y_1=\{i\mid 1\leq i\leq 20\}$,
$\overline y_2=\{1,2,\ldots,12,21,\ldots,28\}$.
Since $\lambda^{(1)}_1=4$, we may assume
$1\in \overline y_j$ for $j=3,4$ and
$1\not\in \overline y_j$ for $j=5,6$. 
Then we must have $\{1,13\}\subset  
\overline y_j$ for $j=1,2,3,4$. 
This impossible.\\
$\bullet $ $30(25)$: $r_2=27$ and $N_{r_2}=10$.
We have $\lambda^{(1)}_1 = 14,
~ \lambda^{(2)}_1=9$ and $\lambda^{(1)}_2(u)=7,
\lambda^{(2)}_2(u)= 9$ for any $u\in X_2$.
Let $Y_{r_2}=\{y_1,\ldots,y_{10}\}$.
Then $r_2=27$ and $\alpha_2=6$,
implies $|\overline y_i\cap \overline y_j|=24$ for any 
distinct $y_i,y_j\in Y_{r_2}$. 
We may assume
$\overline y_1=\{i\mid 1\leq i\leq 27\}$,
$\overline y_2=\{1,2,\ldots,24,28,29,30\}$.
Since $\lambda^{(2)}_1=9$, 
we may assume
$1\in \overline y_j$ for $3\leq j\leq 9$ and
$1\not\in \overline y_{10}$. 
Then $\lambda^{(2)}_2(\{1,30\})=9$ implies 
$\{1,30\}\subset \overline y_j$ for $j=1,\ldots,9$.
But this is impossible since $30\not\in\overline y_1 $. 
\hfill\qed\\

\noindent
{\bf Acknowledgment:} Eiichi Bannai was supported in part by 
NSFC grant No. 11271257. Hideo Bannai was supported in part by
Kakenhi No. 22680013 and No. 25280086.
\\
The authors thank Professor Woodall for kindly answering some
questions by the authors on his paper in 1970.\\

 Eiichi Bannai: Department of Mathematics, Shanghai Jiao Tong University,

800 Dongchuan Road, Shanghai, 200240, China

e-mail: bannai@sjtu.edu.cn
\\
 
 Etsuko Bannai: Misakigaoka 2-8-21, Itoshima-shi, Fukuoka, 819-1136, Japan

e-mail: et-ban@rc4.so-net.ne.jp 
 \\
 
 Hideo Bannai:
 
Department of Informatics, Kyushu University

744 Moto-oka, Nishi-ku, Fukuoka 819-0395, Japan.

e-mail: bannai@inf.kyushu-u.ac.jp

\end{document}